\documentclass[10pt,journal]{IEEEtran}
\usepackage{latexsym, amssymb,amsmath, amsbsy,amsopn, amstext}
\usepackage{latexsym, amssymb,amsfonts,dsfont}
\usepackage{graphicx}\usepackage{subfigure}
\usepackage{tikz}\usepackage{pgf}\usetikzlibrary{arrows,automata}
 \usepackage[latin1]{inputenc}
 \usepackage{verbatim}
\usepackage{amsmath, amsbsy}
\usepackage{amsopn, amstext}
\usepackage{extarrows}
\usepackage{hyperref}
\usepackage{bm}
\usepackage{multicol}
\usepackage{multirow}
\usepackage{tabularx}
\usepackage{cancel, color}
\usepackage{cite}
\usepackage{caption}
\usepackage{float}
\usepackage{stfloats}
\usepackage{ifpdf}
\newtheorem{thm}{Theorem}[section]
\newtheorem{lem}[thm]{Lemma}
\newtheorem{cor}[thm]{Corollary}

\newtheorem{prp}[thm]{Proposition}
\newtheorem{exa}[thm]{Example}
\newtheorem{rem}[thm]{Remark}

\usepackage[thinlines]{easytable}

\allowdisplaybreaks[4]

\def\deq{\stackrel{\Delta}{=}}

\title{ Distributed Stochastic Approximation Algorithm With Expanding Truncations
  }

\author{Jinlong ~Lei, and ~ Han-Fu Chen
\thanks{Jinlong Lei is with the Department of Industrial and Manufacturing Engineering, Pennsylvania
		State University, University Park, PA 16802, USA  (e-mail: jxl800@psu.edu, leijinlong11@mails.ucas.ac.cn).

Han-Fu Chen  is  with the Key Laboratory of Systems and Control, Institute of Systems Science,
Academy of Mathematics and Systems Science,  Chinese Academy of Sciences, Beijing 100190, P. R. China.
(e-mail:  hfchen@iss.ac.cn).}
\thanks{This work was supported by  the National Center for Mathematics and Interdisciplinary Sciences, Chinese Academy of Sciences.}}
\begin{document}
\date{}
\maketitle
\begin{abstract}  In this work, a novel  distributed stochastic approximation  algorithm (DSAA) is proposed  to seek   roots of  the  sum of  local functions, each of which is associated with an agent from the  multiple agents connected over a network. At  each iteration,  each agent updates its estimate for the root utilizing     the noisy observations of its local function and the    information derived from the neighboring agents.  The key difference of  the proposed  algorithm  from the existing ones consists in  the expanding truncations (so it is called as DSAAWET), by which the  boundedness of the estimates can be guaranteed without imposing the growth rate constraints on the local functions.  The  estimates  generated by DSAAWET are shown to
converge  almost surely (a.s.)    to  a consensus  set, which belongs  to a connected subset of  the root set of the sum function.
 In    comparison  with the existing results, we impose weaker conditions on the local functions and on the observation noise.
 We then apply the proposed   algorithm  to two applications,  one from signal processing and the other one from distributed optimization, and
 establish the almost sure convergence.  Numerical  simulation results are also included.
\end{abstract}
\IEEEpeerreviewmaketitle


\section{Introduction }
Distributed algorithms have been  extensively investigated in connection with the problems arising from sensor networks and networked systems in recent years, such as consensus problem  \cite{ren1,murray0},   distributed   estimation \cite{Zhang1,filter1},  sensor localization \cite{localization1},   and distributed optimization \cite{optimization0, nedic,optimization1}.
The distributed algorithms work in the situation,  where the  goal is cooperatively  accomplished by a multi-agent network  with computation and communication abilities allocated in a distributed environment. Their advantages  over the  centralized  approaches for networked problems lie in  enhancing the robustness of the networks, preserving  privacy, and reducing the communication and computation costs.

Stochastic approximation (SA)  was firstly considered by  Robbins and Monro  in 1950s  \cite{RM}   for finding roots of a function with noisy observations, now it is known as the RM algorithm.   Then SA was used by Kiefer and Wolfowitz  \cite{KW} to  estimate the  maximum of an expectation-valued       function  only with noisy function observations.The SA-based schemes have found wide applications in signal processing, communications and adaptive control, see, e.g.,\cite{Kushner,Chen_2002,Chen_Zhao}.
Recently,  many distributed problems are also solved by   SA-based distributed algorithms, e.g.,   distributed  parameter  estimation  \cite{Zhang1,lei2015distributed},  distributed convex optimization over random networks \cite{lei2018asymptotic}, and distributed  non-convex stochastic optimization \cite{Bianchi1}.
  As such, the distributed  variants of SA  recently have   drawn much attention from researchers.

 Distributed stochastic approximation  algorithms  (DSAA)  were  proposed in  \cite{DSA,DSA2}
 to  cooperatively find roots of a function,  being a  sum of  local  functions associated with agents   connected  over a network.  Each agent updates its estimate for the root based on its local information composed of
 the observations of its local function possibly corrupted by noises,  and  
 the  information obtained from  its neighbors.  The  weak convergence   for DSAA with a constant step-size is investigated in \cite{DSA2}, while  the a.s. convergence for   DSAA with a decreasing step-size is  studied     in \cite{DSA}.
Besides, the performance gap between the  distributed and the centralized stochastic approximation algorithms is investigated
in  \cite{Bianchi2}.   The asynchronous and distributed stochastic  approximation algorithms are also addressed  in \cite{Bokar1, Bokar2}  with  the components   separately estimated at different  processors.
However, is noticed that almost all aforementioned  SA-based distributed algorithms   require  rather restrictive  conditions   to guarantee  convergence. For example,  in \cite{DSA}   it is required that each local function is globally Lipschitz continuous and the observation noise is a   martingale difference sequence (MDS)\footnote{Consider an adapted sequence $\{X_t,\mathcal{F}_t\}_{-\infty}^{\infty}$  on a probability space $(\Omega, \mathcal{F}, \mathbb{P})$.
$X_t$ is an MDS if it satisfies $\mathbb{E}[X_t]<\infty$ and $\mathbb{E}[X_t| \mathcal{F}_{t-1}]=0, a.s.,$ for all t.}. However,  these conditions may not hold for some problems, e.g.,  the distributed principal component analysis    and  the distributed gradient-free optimization  to be discussed in Section II.B.  This work aims at   resolving the distributed root-seeking problem  under weaker conditions compared  with those used in \cite{DSA}.

Contributions of this work are as follows. 1)
 We propose  a novel  DSAA with expanding truncations (DSAAWET) over networks with deterministic switching topologies.
 The key difference from the existing algorithms (cf.\cite{DSA,DSA2}) is  the network  expanding truncation mechanism
that  adaptively  defines enclosing  bounds as follows:   the initial value  lies within  the  initial bound;
when   the   estimate  crosses  the   enclosing  bound, it  is reinitialized to the initial point   and     the bound is enlarged.
This together with  the decreasing step-sizes makes as if    the algorithm returns  with   the same initial value but with  smaller step-sizes and larger   bounds.     Then it is shown that   the reinitialization  ceases in a finite number of steps and   the estimates are  bounded.
Further, estimates generated by DSAAWET at all agents converge a.s. to a consensus set,   which is contained in the root set of the sum function. Compared with \cite{DSA},  we neither assume the observation noise to be an MDS nor impose any growth rate constraint and Lipschitz continuity on the local functions.
2) The proposed algorithm is then applied to the distributed principal component analysis    and
to the distributed gradient-free optimization. Their a.s. convergence  is   established as well, while
the algorithms    in  \cite{DSA,DSA2}  might not   be  applicable to  these problems.
 The proposed method  has also been   applied  in \cite{liu2017distributed} and
 \cite{feng2019output}   to the distributed   blind channel identification and   the output consensus of networked nonlinear systems.

  The notations used in this work are  listed in  Table \ref{TAB-nottaion} and the rest of this work is arranged as follows.
The distributed root-seeking problem  is formulated  in Section II along with two motivation examples. DSAAWET is  defined   in Section III and the corresponding  convergent results are   presented as well.  The proof  of   the  main results  is  given in Section IV with some details placed in Appendices.
The proposed algorithm is applied to solve the two application problems  in Section V with    numerical examples included. Some concluding remarks are given in Section  VI.

     \begin{table} [!htb]
\newcommand{\tabincell}[2]{\begin{tabular}{@{}#1@{}}#2\end{tabular}}
\scriptsize
\centering
\begin{tabular}{lcc}
  \hline
  Symbol  &  Definition \\
    \hline \hline
$\| v \| , \| A\|$      & Euclidean ( $l^2$) norm of vector $v$, matrix $A$ \\
  \hline
  $A \geq 0$  &      Each entry of matrix $A$ is nonnegative, \\&and $A$  is called the  nonnegative matrix.    \\
  \hline
  $ \mathbf{I}_m $ &  $m \times m$ identity matrix     \\
    \hline
    $ \mathbf{1} , \mathbf{0} $ &      Vector or matrix   with all entries equal to 1,  0 \\
    \hline
  $X^T$, $X^{-1}$   &  Transpose  of matrix $X$, inverse  of matrix $X$    \\
    \hline
$ col  \{x_1, \cdots, x_m \} $&  \tabincell{c} {   $ col  \{x_1, \cdots, x_m \} \triangleq( x_1^T,\cdots, x_m^T)^T  $\\
 stacking the  vectors or matrices  $x_1,\cdots, x_m$} \\
  \hline
$I_{[\textrm{Inequality}]}  $                        &  Indicator function meaning   it equals 1 if the inequality  \\
  ~  & holds and 0 if the inequality does not hold  \\
  \hline
    $\otimes $        & Kronecker  product  \\
    \hline
     $  a \wedge  b  $   &  min $\{ a, b \}$\\
    \hline
    $J$   &   $  J  \triangleq \{x\in \mathbb{R}^l : f(x)=0 \}$ denotes the root set of $f(\cdot)$ \\ \hline
     $d(x,J)$ &  $ \min_y\{\parallel x-y\parallel: y\in J \} $ \\
     \hline
     $ E [\cdot]$ &  Expectation  operator  \\
     \hline
   &  $ D_{\bot} \triangleq  ( \mathbf{I}_N -\frac{\mathbf{1}\mathbf{1}^T}{N}) \otimes \mathbf{I}_l$,   where  \\
       $D_{\bot} $ & $N$ denotes  the number of agents in the network and
     \\  ~&     $l$ denotes   the  dimension of the root \\ \hline
 $ \delta(t)$ &       $ \delta(t ) \rightarrow 0 \textrm{ as }t \rightarrow 0$.\\ \hline
     $m(k,T)$  & $m(k,T) \triangleq \max \{ m: \sum_{i=k}^m \gamma_i \leq T  \}  $ \\& is an integer-valued function for $T>0$ and integer $k$\\
     \hline
   $\sigma_{i,k}$    & The truncation number of agent $i$ at time $k$ \\ \hline
       $\hat{\sigma}_{i,k}$    &
   $ \hat{\sigma}_{i,k} \deq \max_{  j \in N_{i }(k)} \sigma_{j,k} $, where   $ N_{i}(k) $ is   \\
   &   the set of   neighboring agents of agent  $i$  at time $k$ \\ \hline        $\sigma_{k}$    & The largest truncation number among all agents at time $k$, \\& i.e.,
   $\sigma_{k}  =\max\limits_{i \in \mathcal{V}}\sigma_{i,k}=\max\limits_{i \in \mathcal{V}} \hat{ \sigma}_{i,k}  .$ \\
      \hline
      $\tau_{i,m}$  &  $\tau_{i,m}\deq\textrm{inf}  \{k: \sigma_{i,k}=m \}$   \\
  \hline
 $\tau_m$ &   $\tau_m\deq \min\limits_{i \in \mathcal{V}}\tau_{i,m},$ the smallest time  when at least  \\& one of agents
 has its truncation number reached  $  m$
 \\
 \hline
$\tilde{\tau}_{j,m}$ &   $\tilde{\tau}_{j,m}\deq \tau_{j,m} \wedge \tau_{m+1}$ \\
  \hline
   \end{tabular}
\caption{Notations}\label{TAB-nottaion}
\end{table}

 \section{ Problem Formulation and Motivations}

In this section, we first formulate the distributed root-seeking problem with the related communication model.
Then   we   give two motivation  problems that cannot be solved by the existing algorithms,
but can   be solved by DSAAWET   to be proposed in this work.

\subsection{  Problem Statement}
Consider the case where all agents in a network  collectively search the root of the sum function given by
\begin{equation}\label{summ1}
\begin{array}{lll}
    f(\cdot)= \frac{1}{N}\sum_{i=1}^N f_i(\cdot),
    \end{array}
   \end{equation}
    where  $f_i(\cdot):\mathbb{R}^l \rightarrow  \mathbb{R}^l$  is the local function assigned to agent $i$ and can only be observed by agent $i.$       Let   $J \triangleq  \{x\in \mathbb{R}^l : f(x)=0 \}$ denote the root set of $f(\cdot)$.

Time is slotted at $k=0,1,2,\dots$
  For any $   i \in \mathcal{V}$,  denote by $x_{i,k} \in \mathbb{R}^l$
   the  estimate for the   root of   $f(\cdot) $ given by  agent  $i  $
   at time $k$.  Agent $i$  at time $k+1$ has its local  noisy observation
     \begin{equation}\label{ob}
    O_{i,k+1}= f_i(x_{i,k})+\varepsilon_{i,k+1},
     \end{equation}
 where  $\varepsilon_{i,k+1}$  is the observation noise.
Agent $i$ is required to update  its
 estimate   $ x_{i,k}$  on the basis of   its local observation and
the  information  obtained from its    neighbors.

The  information exchange  among  the  $N$ agents at time $k$ is described by  a  digraph
   $\mathcal{G }(k) = \{ \mathcal{V }, \mathcal{E }(k)\}$, where  $\mathcal{V }=\{ 1,\cdots,N\}$ is the node set with  node $i$ representing agent $i$;
 $\mathcal{E }(k) \subset  \mathcal{V } \times \mathcal{V } $ is the edge set with    $(j,i)\in\mathcal{E }(k)$ if   agent  $i$ can  obtain  information from agent $j$ at time $k$.  Let the associated  adjacency matrix
 be  denoted by $W(k)  =[ \omega_{ij}(k)]_{i,j =1}^N$, where  $ \omega_{ij}(k)>0$  if  and only if  $(j,i )\in \mathcal{E }(k)$,   and  $ \omega_{ij}(k)=0$, otherwise. Denote  by $ N_{i}(k) \triangleq \{ j \in \mathcal{V}:(j,i )\in \mathcal{E }(k) \}$  the set of neighboring agents of agent  $i$  at time $k$.

A  time-independent   digraph  $\mathcal{G}=\{ \mathcal{V}, \mathcal{E}\}$  is  called  strongly connected
if  for any $  i,j\in \mathcal{V}$ there exists a directed path from  $i$ to $j$.    By this  we mean
 a  sequence  of edges   $ (i,i_1),(i_1,i_2),\cdots,(i_{p-1},j)$  in the digraph
 with  distinct  nodes $ i_m \in \mathcal{V}~~\forall m: 0 \leq m \leq p-1$,  where $p$ is called  the length of the directed path.
A nonnegative square   matrix $A $   is  called doubly stochastic
 if  $A\mathbf{1} =\mathbf{1}$ and  $\mathbf{1}^T A =\mathbf{1}^T$.

    \subsection{   Motivation  Examples}

We now  give two motivation examples   that  cannot be solved by the existing distributed stochastic approximation algorithm \cite{DSA}. We will return to these examples in Section V
to show  that  DSAAWET  to be proposed in this work   can  solve them.

\subsubsection{Distributed Principal Component Analysis}
In  signal processing and pattern  recognition,    effective  clustering of  large data sets is an important  objective.
 Principal component analysis  (PCA) is a powerful technique  to process  the multivariate data by  constructing a concise data representation through   computing the dominant eigenvalues  and the corresponding eigenvectors of the data covariance matrix \cite{Oja}.
  For recent years,  the  distributed computation of  PCA has attracted much attention  from  researchers  (cf. \cite{DPCA1,macua2010consensus,DPCA2,korada2011gossip,morral2012asynchronous,meng2012distributed}).
Roughly speaking,  there are two different types of distributed PCA:  (i)
Samples composed of all features    collected  from the wireless  sensor networks (WSNs) \cite{DPCA1,macua2010consensus}.  (ii)  Agents  in the network   collect    observations for    different  features  of the whole features \cite{DPCA2,korada2011gossip,morral2012asynchronous,meng2012distributed}.
The global PCA is computed in   \cite{DPCA1}  for distributed updating data sets  without convergence analysis, while     the   distributed PCA of a covariance matrix of the  fixed data sets  is  addressed in  \cite{macua2010consensus}.  The authors of  \cite{DPCA2}  concentrate  on computing  principal components   of  a covariance matrix  observed with  noise,  propose an  algorithm with each sensor estimating   the corresponding entry of the principal eigenvector, and  show  that the  associated ODE is close to that of Oja's  \cite{Oja}.  A distributed PCA algorithm  is proposed  on the directed Gaussian
graphical models in \cite{meng2012distributed}.
The goal of  \cite{korada2011gossip,morral2012asynchronous}  is to distributively compute
   the principal components   of   a symmetric weighted adjacency matrix  associated with the connected network.

We consider the distributed setting similar to that  considered  in \cite{DPCA1,macua2010consensus}
 aiming at proposing  a distributed algorithm for the global PCA and    establishing  its a.s. convergence.
The global data matrix $u_k$ at time $k$ is distributed among  $N$ sensors  that are spatially distributed  in the network:
$$u_k= \begin{pmatrix}
    u_{1,k}  \\
       \cdots \\
      u_{N,k}
\end{pmatrix} ,$$
where  $u_{i,k} \in \mathbb{R}^{p_i\times l}$ is collected by sensor $i.$  The rows  of $u_k$ denote the observations
while the  columns denote the features.
    Assume that for each $i \in \mathcal{V}$, $\{u_{i,k}\}$ is an i.i.d. sequence with zero mean.
   Denote by  $A \triangleq  E[u_k^Tu_k]  \in \mathbb{R}^{l\times l}$    the   covariance  matrix    of  $u_k$.
The primary objective is to  estimate  the unit eigenvector    corresponding to    the largest eigenvalue of $A$.
Because the unit eigenvector corresponding to the largest eigenvalue of $A$
 belongs to the set of the  nonzero roots  of  the function  $  f(x)=Ax-(x^T Ax)x$ on the unit sphere,   then the   distributed   PCA might be  solved     by virtue of finding   roots of  the sum function defined as follows:
  $$f(x)= \sum_{i=1}^N f_i(x) ,$$
  where $ f_i(x)=A_ix-(x^T  A_ix)x$ with  $ A_i\triangleq E[u_{i,k}^Tu_{i,k}]$.

    Denote  by  $x_{i,k}$   the  estimate given by agent $i$ at time $k$ for the unit eigenvector  corresponding to
     the largest eigenvalue of matrix $A$. Since $A_i$ cannot be directly derived,
      by replacing   $A_i$ with its sample  data $A_{i,k}=u_{i,k}^Tu_{i,k}$ at time $k$,
      the local observation of  agent $i$ at  time $k+1$ is defined  as follows:
  \begin{equation}\label{pca}
    O_{i,k+1}\triangleq A_{i,k}x_{i,k}-(x_{i,k}^T A_{i,k}x_{i,k})x_{i,k} .
  \end{equation}
    Then the distributed principal component analysis  is in the distributed root-seeking form.
  Noting that   $f_i(\cdot)$ contains a cubic term,    the local function is not globally Lipschitz continuous.
 Though  the existing distributed stochastic approximation algorithms  might be still applicable since   eigenvectors can be chosen to have unit-norm by projecting the estimates onto the unit hypersphere,  the exact convergence to the true eigenvector for the distributed data set with updating has not been fully addressed yet.  In Section V, we will show that    DSAAWET     can solve this problem and ensure the a.s. convergence.

   \newcounter{mytempeqncnt}
\begin{figure*}[!b]
\normalsize
\setcounter{mytempeqncnt}{\value{equation}}
\setcounter{equation}{5}
\hrulefill
 \begin{align}
  & \sigma_{i,0}=0,~~~     \hat{\sigma}_{i,k} \deq \max_{  j \in N_{i }(k)} \sigma_{j,k},  \label{comp1} \\
    &  x'_{i,k+1}= \Big( \sum_{ j \in  N_{i }(k)}  \omega_{ij}(k) ( x_{j,k} I_{[ \sigma_{j,k} = \hat{ \sigma}_{i,k}]}+
      x^* I_{[ \sigma_{j,k} <  \hat{\sigma}_{i,k}]})   + \gamma_{ k } O_{i,k+1} \Big) I_{[\sigma_{i,k}= \hat{ \sigma}_{i,k}]}
 + x^* I_{[\sigma_{i,k}< \hat{ \sigma}_{i,k}]}, \label{comp2} \\
   &   x_{i,k+1}=x^* I_{[   \parallel x'_{i,k+1} \parallel > M_{\hat{\sigma}_{i,k}} ]}
      +x'_{i,k+1}I_{[   \parallel x'_{i,k+1} \parallel \leq M_{\hat{\sigma}_{i,k}}] },  \label{step1}\\
   & \sigma_{i,k+1}= \hat{\sigma}_{i,k}+  I_{\left[ \parallel x'_{i,k+1} \parallel > M_{\hat{\sigma}_{i,k}}\right]},  \label{step2}
\end{align}
\setcounter{equation}{\value{mytempeqncnt}}
\end{figure*}

\subsubsection{Distributed Gradient-free  Optimization}

Consider a multi-agent network of $N$ agents,  for which  the objective   is
  to  cooperatively solve the following  optimization problem:
  \begin{equation}\label{opt1}
     \min_x \quad c(x)=\sum_{i=1}^N c_i(x),
  \end{equation}
  where  $c_i(\cdot) : \mathbb{R}^l \rightarrow \mathbb{R}$ is the local objective function of  agent   $i$,
  and $c_i(\cdot)$ is only known by   $i$ itself.     Assume the optimization problem \eqref{opt1} has solutions.   Consider the  case  where   gradients of the cost functions are unavailable but   the cost functions   can be observed with noises.
Then  the finite time difference of the cost functions  can be adopted to estimate the gradient, see \cite{KW,Chen_1999,gradienfree2}.  This problem is referred    to as the  gradient-free optimization   \cite{gradienfree1}.
 The gradient-free methods   for nonsmooth distributed  optimization  with  convex set constraint  are considered  in  \cite{Yuan_NN_15} \cite{li2015gradient},  where it is assumed that each local objective  function  is convex and   Lipschitz continuous over the convex set.
The underlying network topology is modeled as time varying in  \cite{Yuan_NN_15} while  as fixed in \cite{li2015gradient}.
In both  \cite{Yuan_NN_15} and \cite{li2015gradient},   the Gaussian smoothing technique is used,  and the convergence   to an approximate solution   within the error level depending on the smoothing parameter and the Lipschitz constant of local objective function is established.    In contrast, we consider the  distributed  smooth  optimization  with the aim of  designing a  distributed gradient-free algorithm that can achieve  a.s. convergence while without  requiring each objective function   be convex and Lipschitz  continuous.

 Assume that  the cost functions  $c_i(\cdot),i=1,\cdots N$
  are differentiable and the global function $c(\cdot)$ is convex.  Then solving  the problem \eqref{opt1}  is equivalent to  finding  roots of  the function $ f(x)= \sum_{i=1}^N f_i(x)/N$ with $f_i(\cdot)\triangleq- \nabla c_i(\cdot)$,
    where $\nabla c_i(\cdot)$ denotes  the gradient function of $c_i(\cdot)$.
Here the  randomized KW method proposed in  \cite{Chen_1999} is  adopted to estimate   gradients
of the  cost functions.  Denote  by $x_{i,k}$  the   estimate of  the solution to the
 problem  \eqref{opt1} given by agent $i$  at time $k$.   Let $\Delta_{i,k} \in \mathbb{R}^l,k=1,2,\dots$ be
   a sequence of mutually independent  random vectors with  each component independently
    taking values $\pm1$ with   probability $\frac{1}{2}$.
     The   observation  of  function $- \nabla c_i(\cdot)$ at point $x_{i,k}$ is  constructed  as  
\begin{equation}\label{opt2}
   O_{i,k+1}=-\frac{c_i(x_{i,k }^{+})+\xi_{i,k+1}^{+}-c_i(x_{i,k }^{-})-\xi_{i,k+1}^{-}}{2\alpha_k} \Delta_{i,k},
\end{equation}
   where $x_{i,k }^{+}=x_{i,k}+\alpha_k\Delta_{i,k}$, $x_{i,k }^{-}=x_{i,k}-\alpha_k\Delta_{i,k}$, $\alpha_k >0  $ for any $k \geq 0$,    $\xi_{i,k+1}^{+}$  and $\xi_{i,k+1}^{-}$ are the observation noises of  the cost function
   $c_i(\cdot)$ at   points  $x_{i,k }^{+}$ and $x_{i,k }^{-}$, respectively.

 Thus,  we have transformed the distributed gradient-free optimization   to the  distributed root-seeking form.
 It is  shown in \cite{Chen_1999} that   for any $i \in \mathcal{V}$, the observation noise  $ \varepsilon_{i,k+1}=O_{i,k+1}+\nabla c_i(x_{i,k}) $ is  not an  MDS, and hence DSAA  proposed in \cite{DSA}  cannot be directly    applied to this problem.

  The problems arising from   this kind of  examples motivate us to propose DSAAWET, by which the two motivation examples can be resolved as to be shown in Section V.

\section{DSAAWET and Its Convergence}
In this section, we  define DSAAWET and formulate the main results of this work.

 \subsection{ DSAAWET}

Let us first explain the idea of   expanding truncations.  In many important problems, the sequence of estimates generated  by the RM algorithm may not be bounded, and it is hard to  define in advance  a region  the sought-for parameter belongs to.  This motives us to adaptively define truncation bounds as follows.
 When the   estimate  crosses  the  current  truncation bound,  the estimate  is reinitialized to  some pre-specified point  $x^*$ and   at the same  time  the  truncation bound is enlarged.  Due to the decreasing step-sizes,   this operation makes as if we  rerun the algorithm  with initial value $x^*$,   with smaller step-sizes and larger  truncation bounds.     The  expanding truncation mechanism  incorporating with some verifiable conditions makes the reinitialization  cease in a finite number of steps, and hence makes the estimates bounded.
As results,  after a finite number of steps the algorithm runs as the RM algorithm.   This  is well explained in
 \cite{Chen_2002} and    the related references therein.

We  apply the idea of expanding truncation  to the distributed estimation, i.e., for the case $N>1$.  Let $x_{i,k}$  denote agent $i$' estimate for the root of $f(\cdot)$ at time $k,$  and $\sigma_{i,k}$  denote  the number of truncations for  agent $i$  up-to-time $k$.  We now  define DSAAWET  (Eqns. \eqref{comp1}-\eqref{step2}) with  any initial values $x_{i,0},$   where   $ O_{i,k+1}$  is defined by \eqref{ob}, $\{ \gamma_k\} $  is the   steplength,
  $x^*  \in  \mathbb{R}^l$ is a fixed     vector known to all agents,  
    $\{M_k\}_{k \geq 0} $   is  a  sequence   of positive numbers increasingly diverging to infinity with     $M_0 \geq \| x^* \|$,  and $  M_{\hat{\sigma}_{i,k}}$
serves as the  the truncation bound  when the $(k+1)$th estimate  for agent  $i$ is generated.
 Note that during an iteration, each agent  uses  its neighbors'  information including estimates for
  the root and the truncation number.  

The algorithm  \eqref{comp1}-\eqref{step2}  is performed  according to the following three steps.

\noindent  (i)  \textbf{Max consensus of  truncation numbers} (Eqn. \eqref{comp1}).  At time $k$, the update of  agent $i$
may or may not be truncated. This yields that at time $k,$ $\sigma_{i,k}$ may be different for $i=1,\cdots, N$.
  We then  set the truncation number $\hat{\sigma}_{i,k}$ to be the largest one among $\{ \sigma_{j,k},~j \in \mathcal{N}_i(k) \} $ of  its neighboring agents     as indicated by \eqref{comp1}. It is worth noting that we do not require the truncation number  of  each agent  to reach a  max consensus in finite time. Nevertheless, we   will  prove  in Lemmas   \ref{lemma4} and \ref{lem001}   that  all agents' truncation numbers will asymptotically  reach a consensus to some positive integer. 
   
 \noindent   (ii) \textbf{ Average consensus + innovation update}  (Eqn. \eqref{comp2}).  At time $k+1$,  agent  $i$ produces an  intermediate  value  $x'_{i,k+1}$   by \eqref{comp2}.
 If  $\sigma_{i,k}< \hat{ \sigma}_{i,k} $, then  $x'_{i,k+1}=x^{*}$.
 Otherwise,  $x'_{i,k+1}$    is a combination of the consensus part and the innovation part,
 where the consensus part is  a weighted average of estimates derived at its  neighbors,
  and the innovation part processes its local current observation.

\noindent  (iii)   \textbf{Local truncation judgement}  (Eqns. \eqref{step1}-\eqref{step2}).  If  $  x'_{i,k+1}  $ remains inside its local
 truncation bound $M_{\hat{\sigma}_{i,k}}$, i.e., $\parallel x'_{i,k+1} \parallel  \leq  M_{\hat{\sigma}_{i,k}} $, then $ x_{i,k+1}= x'_{i,k+1}$ and $\sigma_{i,k+1}= \hat{\sigma}_{i,k}$.  If
       $x'_{i,k+1}$  exits from the  sphere with radius $M_{\hat{\sigma}_{i,k}}$,  i.e., $\parallel x'_{i,k+1} \parallel >M_{\hat{\sigma}_{i,k}} $, then    $ x_{i,k+1} $ is pulled back   to the pre-specified   point $ x^*$,    meanwhile,  the truncation number  $\sigma_{i,k+1}$ is increased to  $  \hat{\sigma}_{i,k}+1$.   

Denote  by
\setcounter{equation}{9}
    \begin{equation}\label{largest}
\sigma_{k}   \deq \max\limits_{i \in \mathcal{V}}\sigma_{i,k},
\end{equation}
 the largest truncation number among all agents at time $k$. If $N=1,$ then  by denoting
  $   \hat{ \sigma}_{i,k} =\sigma_{i,k}\triangleq \sigma_k,~ x_{i,k} \triangleq x_k, ~x_{i,k}'\triangleq x_k', ~O_{i,k}\triangleq O_k,$  it is  seen that  \eqref{comp1}-\eqref{step2} becomes
   \begin{equation}
   \begin{split}
    &  x'_{k+1}= x_k+ \gamma_{ k } O_{k+1}  , \\
   &   x_{k+1}=x^* I_{\left[   \parallel x'_{k+1} \parallel > M_{ \sigma_k }\right]}
      +x'_{k+1}I_{\left[   \parallel x'_{k+1} \parallel \leq M_{ \sigma_k}\right]},   \\
   & \sigma_{k+1}= \sigma_k+  I_{\left[ \parallel x'_{k+1} \parallel > M_{ \sigma_k}\right].} \nonumber   \end{split}
\end{equation}
The above algorithm  is called the stochastic approximation algorithm with expanding truncations (SAAWET), which requires possibly the weakest conditions for its convergence among various modifications of the RM algorithm  \cite{Kushner, Chen_2002, Ben, Bokar}. Thus, the advantages   of SAAWET over the  RM  algorithm might remain for its distributed variant in the case $N>1.$

\begin{rem}\label{r1}
It is noticed that $\sigma_{i,k+1}\geq \hat{\sigma}_{i,k} \geq \sigma_{i,k} ~~\forall k \geq 0$  by \eqref{comp1} and \eqref{step2}. Further,  it is concluded that
\begin{equation}\label{sub1}
x_{i,k+1}= x^*\textrm{  if  }\sigma_{i,k+1} > \sigma_{i,k}.
\end{equation}
This can be seen from the following consideration:
i)  If    $ \sigma_{j,k} \leq  \sigma_{i,k} ~\forall j \in N_{i} (k)$, then
  from \eqref{comp1}  we derive  $ \hat{\sigma}_{i,k} =   \sigma_{i,k}$.
 Since   $ \sigma_{i,k+1} >  \sigma_{i,k}$, by \eqref{step2} it follows that
      $\parallel x'_{i,k+1} \parallel > M_{\hat{\sigma}_{i,k}}$,
and hence from    \eqref{step1}  we derive $x_{i,k+1}=x^*.$
ii) If there exists $  j \in N_i(k)$ such that  $ \sigma_{j,k} >  \sigma_{i,k}$,
 then from  \eqref{comp1}  we derive $ \hat{\sigma}_{i,k} =   \max_{  j \in N_{i }(k)} \sigma_{j,k} > \sigma_{i,k},$
 and  from \eqref{comp2} we have  $x'_{i,k+1}=x^*$.
 Consequently, by  \eqref{step1}  we have  $x_{i,k+1}=x^*$.
    \end{rem}

\subsection{Assumptions}
 We list the   assumptions to be used.

 A1   $ \gamma_k>0, \gamma_k   \xlongrightarrow [k \rightarrow \infty]{}   0 $, and  $ \sum\limits_{k=1}^\infty \gamma_k=\infty$.

A2 There exists a continuously differentiable function  $v(\cdot  ) :
\mathbb{R}^l\rightarrow \mathbb{R}$ such that
\begin{flalign}
& ~~~\textrm{a)} ~~~~~~~~~~~~\sup_{{\delta \leq d(x,J  )\leq\Delta}   }  f^T(x)v_x(x)<0&  \label{4.1}
\end{flalign}
for any  $ \Delta > \delta>0$,           where   $v_x(\cdot)$
denotes the gradient of $v(\cdot) $  and $d(x,J )=  \min_y\{\parallel x-y\parallel: y\in J \} $,

b)  $v(J )\triangleq \{v(x): x\in J  \}$ is nowhere dense,

c) $  \| x^* \| < c_0  \textrm{ and } v(x^*) < \textrm{inf }_{\| x \|=c_0} v(x) $ for some positive constant
$c_0 $,   where $x^*$ is used in   \eqref{comp2} \eqref{step1}.

A3   The local functions  $f_i(\cdot)~~ \forall i \in \mathcal{V}$ are continuous.

Let us  explain    the conditions.  A1 is a standard assumption for stochastic approximation, see \cite{Chen_2002,Bokar}.
 A2  a)  means that $v(\cdot)$ serves as a Lyapunov function for the differential equation $\dot{x}=f(x).$
It is noticed that  A2  b) holds if $J$ is finite,   and A2 c)  takes place  if $v(\cdot)$ is radially  unbounded.

 The following condition  A4 imposed on  the communication graphs is taken from  \cite{optimization0},
  which the readers are referred to for detailed explanations.

A4 a)   For any $k\geq 0,$ $W(k) $ is a doubly stochastic matrix;

~~~ b) There exists a constant  $0< \eta <1$  such that
 $$ \omega_{ij}(k) \geq \eta   \quad  \forall j \in  \mathcal{N}_i(k) ~~\forall i \in \mathcal{V}~~\forall k \geq 0;$$

 ~~~ c)  The digraph  $\mathcal{G}_{\infty}=\{ \mathcal{V}, \mathcal{E}_{\infty}\}$
 is strongly connected,
 where $\mathcal{E}_{ \infty}$ is the set of edges $(j,i)$ such that $j$
is a neighbor of $i$ which  communicates with $i$ infinitely often, i.e.,
$$\mathcal{E}_{ \infty} =
 \{ (j,i):  (j,i) \in \mathcal{E}(k)  \textrm{ for infinitely   many indices~} k \} ;$$

~~~  d) There exists  a positive integer $B $  such that for every $(j,i) \in {\cal E}_\infty$ ,
agent $j$ sends  the information to the neighbor  $j$ at
least once every $B$ consecutive time slots, i.e.,
 \begin{equation*}
    (j,i) \in \mathcal{E} (k) \cup \mathcal{E}(k+1)   \cup  \cdots \cup \mathcal{E}(k+B-1)
 \end{equation*}
 for all $(j,i) \in {\cal E}_\infty$ and any $k\geq 0.$

\begin{rem}
 A4 a)  requires the adjacency matrices of the digraphs to be doubly stochastic, which   does not hold for all digraphs.
Nevertheless,   a   necessary and sufficient condition is provided  in \cite{gharesifard2010does}  for a digraph to be doubly stochasticable. A4 b) states that each agent gives significant
weights to its  own and  the neighbors' states. It  is explained  in \cite{optimization0} that A4 c) is equivalent
to the assumption that the composite directed graph $ \{ \mathcal{V}, \cup_{l\geq k} \mathcal{E}_{k} \}$ is  strongly connected for all $k$.  A4 d) implies that  for any   $i,j \in \mathcal{V}$,   the information of   $j$ can be  propagated to
 $i$ in less than $B  (N-1)$ steps.   Define
$\Phi(k,k+1) =\mathbf{I}_N $ and $  \Phi(k,s)= W(k) \cdots W(s)$ for all $   k \geq s. $
Then by   \cite[Proposition 1]{optimization0},  there exist  $c>0$ and $0< \rho<1$  such that \begin{equation}\label{graphmatrix}
  \parallel \Phi(k,s) - \frac{1}{N} \mathbf{1} \mathbf{1}^T \parallel \leq c\rho^{k-s+1} \quad \forall k \geq s .
\end{equation}
\end{rem} 

  Before  introducing conditions on noise,  let us denote by $(\Omega,\mathcal{F},\mathbb{P})$ the probability
space.  For  any  $i \in \mathcal{V}$, let $ \varepsilon_{i,k+1}:(\mathbb{R}^l \times \Omega, \mathbb{B}^l \times \mathcal{F})
\to (\mathbb{R}^l \times \mathbb{B}^l )$ be a measurable function, where $ \mathbb{B}$  is a  Borel $\sigma-$algebra. Let the  noise be given
by  $ \varepsilon_{i,k+1}=\varepsilon_{i,k+1}(x_{i,k}(\omega),\omega),\omega\in\Omega.$
Thus,   the state-dependent noise is considered.
For notational simplicity,  $ \varepsilon_{i,k+1}$ is used to denote $ \varepsilon_{i,k+1}\left(x_{i,k}(\omega),\omega\right)$
 for the sample path $\omega \in \Omega$ under consideration throughout this work.

A5   For the sample path  $\omega$ under consideration,  the following assertions  hold for  any  $i \in \mathcal{V}$:
\begin{flalign}
&~~~\textrm{a)} ~  \gamma_k \varepsilon_{i, k+1}  \xlongrightarrow [k \rightarrow \infty]{} 0,\textrm{  and }& \nonumber \\
 &~~~\textrm{b)} ~\lim_{T \rightarrow 0} \limsup_{k \rightarrow \infty }\frac{1}{T}
    \parallel  \sum_{m=n_k}^{m(n_k,t_k)}\gamma_{m } \varepsilon_{i,m+1} I_{ [\parallel x_{i,m} \parallel \leq K ]} \parallel =0
    \nonumber \\
& \quad~~  ~~~~\forall t_k \in [0,T]\textrm{ for   any sufficiently large } K  &  \nonumber
\end{flalign}
along indices $\{n_k \}$   whenever  $\{x_{i,n_k}\}$ converges, where  $m(k,T) \triangleq  \max
\{ m: \sum_{i=k}^m \gamma_i \leq T  \}  $.

\begin{rem}
   It is noticed that A5 b) is convenient for dealing with state-dependent noise. The indicator function   $I_{ [\parallel x_{i,m} \parallel \leq K ]} $ in the condition     will
     be  casted away if the observation noise does not depend on the estimates.    However, for  the  state-dependent   noise, before establishing the boundedness of $\{ x_{i,k}\}$, the condition with an indicator function included  is easier to be verified. It is worth noting that in A5 b)  we do not assume  existence of a convergent  subsequence of
    $\{x_{i,k}\} $ for any $ i$, we only require A5 b) hold    along  indices of any convergent
    subsequence if  exists.   Verification of A5 b) along
   convergent subsequences is much easier than  that   along the whole sequence.
If  $\{ \varepsilon_{i,k}\}$  can be decomposed into two parts
 $\varepsilon_{i,k}=\varepsilon_{i,k}^{(1)}+\varepsilon_{i,k}^{(2)}$ such that
 $ \sum_{k=0}^{\infty } \gamma_{k } \varepsilon_{i,k+1}^{(1)} I_{ [\parallel x_{i,k} \parallel \leq K ]} < \infty
\textrm{ and  }\varepsilon_{i,k}^{(2)} I_{ [\parallel x_{i,k} \parallel \leq K ]} \xlongrightarrow [k\rightarrow \infty] {} 0,$
then  A5  b)  holds.   So,  A5 holds    when the observation noise is an i.i.d. sequence or an MDS with   bounded  second moments if $ \sum\limits_{k=1}^\infty \gamma_k^2<\infty$.
\end{rem}

\subsection{  Main    Results}
Define the vectors  $ X_k \deq col\{ x_{1,k} ,\cdots, x_{N,k}\}, ~\varepsilon_k\deq   col \{
    \varepsilon_{1,k}  , \cdots,   \varepsilon_{N,k}  \} ,
    F(X_k)\deq col \{ f_1(x_{1,k}) , \cdots,f_N(x_{N,k})\} $.
    Denote by   $ X_{\bot, k}\deq D_{\bot}  X_k$ the disagreement vector of $X_k$
     with $ D_{\bot} \triangleq  ( \mathbf{I}_N -\frac{\mathbf{1}\mathbf{1}^T}{N}) \otimes \mathbf{I}_l$,    and  by   $x_k
=\frac{1}{N} \sum_{i=1}^N x_{i,k}$     the average of the estimates derived at   all agents at time  $k$.

   \begin{thm}\label{thm1}
Let $ \{ x_{i,k}\} $ be produced by \eqref{comp1}-\eqref{step2}  with an arbitrary  initial value  $   x_{i,0} $.
 Assume  A1-A4 hold. Then for    any sample path  $\omega$  where  A5    holds,  we have the following:

i)     $\{x_{i,k}\}$ is bounded and   there exists a positive integer  $k_0$
   possibly depending  on $\omega$ such   that
    \begin{equation} \label{notr}
 x_{i,k+1}=  \sum_{ j \in  N_{i }(k)}  \omega_{ij}(k)  x_{j,k}  + \gamma_{ k } O_{i,k+1}      ~~ \forall k \geq  k_0 ,
  \end{equation}
 or in the following  compact form:
   \begin{equation}\label{centalform}
    X_{k+1}= \left( W (k) \otimes \mathbf{I}_l \right)  X_k+ \gamma_k (F(X_k)+\varepsilon_{k+1})  ~\forall k \geq k_0;
 \end{equation}
 \begin{flalign}
&~~~\textrm{ii)} ~ ~~~~  X_{\bot, k} \xlongrightarrow [k \rightarrow \infty]{} \mathbf{0}~~~\mbox{and}~
   \quad d(  x_k,J  ) \xlongrightarrow [k \rightarrow \infty] {} 0; &\label{res1}
\end{flalign}

iii)   there exists  a connected  subset $J^{*} \subset J$ such that
  \begin{equation}\label{connected}
    d(  x_k,J^{*}  ) \xlongrightarrow [k \rightarrow \infty] {} 0.
  \end{equation}
\end{thm}

The proof of Theorem \ref{thm1} is  presented in Section  IV.  It is noticed from  the proof
that  for deriving  i),  condition A5 a) is not required.
    Theorem \ref{thm1} establishes that  the  sequence  $\{X_k\}$ is  bounded;   the
 algorithm  \eqref{comp1}-\eqref{step2} finally  turns to   a   standard DSAA without truncations; and
    the   estimates  given by all agents   converge a.s. to a consensus set,  contained in a connected subset of  the  root set $J $,    when A5 holds for almost all sample paths for all agents.
As a consequence,  if $J$ is not dense in any connected set, then $x_k$ converges to a  point in $J$.
However,   it is   unclear how  $x_k$ behaves when $J$
is dense in some connected set.    This problem was investigated  for
the centralized  algorithm in \cite{fang}.
 Based on Eqn. \eqref{notr}, we might be able to investigate the  rate of convergence of  DSAAWET similar to the standard DSAA,  while we leave  it for future research  due to  the page limitation.

\begin{rem}   Compared      with \cite{DSA}, we  impose  weaker conditions on the local functions and on  the observation noise. In fact,  we only require the local functions  be continuous, while conditions ST1 and ST2 in   \cite[Theorem 3]{DSA} do not allow the  functions to grow faster than   linearly. In addition,
we do not require the observation noise  to  be an MDS as in  \cite{DSA}.   As shown in \cite{Chen_2002, Chen_Zhao}, A5 is probably the weakest requirement for  the noise since it is also necessary for convergence whenever the root $x^0$ of $f(\cdot)$ is a singleton and $f(\cdot)$ is continuous at $x^0 $. Compared with the random  communication
  graphs used in \cite{DSA}, here we use  the deterministic  switching graphs to
  describe the communication relationships among agents.
\end{rem}

\section{Proof  of Main Results } Prior to  analyzing $\{x_{i,k}\}$, let us recall the convergence analysis
  for  SAAWET, i.e., DSAAWET with $N=1$.  The key step in the analysis is to establish
  the boundedness of the estimates, or to show that truncations cease in a
  finite number of steps.  If the number of truncations increases  unboundedly,
  then SAAWET is pulled back to a   fixed vector $x^{*}$ infinitely many times.
  This produces  convergent subsequences from the estimation  sequence.  Then the condition A5 b)
   is applicable and it incorporating  with A2 yields a contradiction. This proves the boundedness of the estimates.

    Let us try to use this approach to prove the boundedness of
    $x_k=\frac{1}{N} \sum_{i=1}^{N} x_{i,k}$ with $\{x_{i,k}\}$ generated by \eqref{comp1}-\eqref{step2}.
    In   the case     $  \sigma_k \xlongrightarrow [k \rightarrow \infty ]{} \infty$, we  have   $ \lim\limits_{k \rightarrow \infty } \sigma_{i,k} = \infty~~ \forall i \in \mathcal{V}$ by Corollary \ref{cor}
given  below.  Then from Remark \ref{r1} it is known that the estimate $x_{i,k+1}$ given by agent $i$
      is  pulled back to  $x^{*}$  when the truncation occurs at time $k+1$.
      This means  that $\{x_{i,k}\}~\forall i \in \mathcal{V}$  contains    convergent subsequences.
      However,  $\{x_k\}$ may still not  contain any convergent subsequence to make A5 b) applicable. This is
      because truncations may occur at different times for different $i \in \mathcal{V}.$ Therefore, the conventional approach used for convergence analysis of SAAWET cannot directly be applied to the algorithm \eqref{comp1}-\eqref{step2}.

To overcome the difficulty, we first introduce   auxiliary sequences $ \{\tilde{x}_{ i,k}\} $ and
 $ \{\tilde{\varepsilon}_{i,k+1} \} $ for any  $i \in \mathcal{V}.$
It will be shown in Lemma \ref{lemau2} that  $ \{\tilde{x}_{i,k}\} $ satisfies the recursions  \eqref{algorithm01}-\eqref{algorithm03}, for which the truncation bound at time $k$ is the same $M_{\sigma_k}$ for all agents  and  the estimates $\tilde{x}_{i,k+1} ~ \forall i \in \mathcal{V}$ are  pulled back to  $x^{*}$  when  $\sigma_{k+1}>\sigma_k$.
As a result,  the auxiliary sequence $\{ \widetilde{X}_k\}$ has convergent subsequences, where   $ \widetilde{X}_k\deq col \{ \tilde{x}_{1,k}, \cdots,  \tilde{x}_{N,k} \}.$
Besides, it will be shown in  Lemma \ref{lemma1} that the noise   $\{\tilde{\varepsilon}_{k} \}$   satisfies a  condition    similar to A5 b) along   any convergent subsequence of  $\{ \widetilde{X}_k\}$, where   $\tilde{\varepsilon}_{k} \deq col \{ \tilde{\varepsilon}_{1,k}, \cdots, \tilde{\varepsilon}_{N,k} \}$. To borrow  the  analytical  method  from the centralized stochastic  approximation algorithm, we  rewrite the algorithm   \eqref{algorithm01}-\eqref{algorithm03}    in the centralized form   \eqref{centrelized2}  with observation noise $\{ \zeta_{m+1}\}$.     By  the results given in Lemma \ref{lemau2}  and  Lemma \ref{lemma1}, it is shown  in Lemma  \ref{lemnoise}   that   the noise sequence  $\{ \zeta_{m+1}\}$ satisfies  \eqref{err2}  along   convergent   subsequences, which is  similar to A5 b) when $N=1$.   Then by algorithm  \eqref{centrelized2} and the noise property  \eqref{err2},  we   show that the  numbers  of truncations for all agents converge to the same finite value, and that $\{\tilde{x}_{i,k}\}_{i \in \mathcal{V}}$  reach a consensus located in  the root set.
Thus,   $\{x_{i,k}\}$ and $\{\tilde{x}_{i,k}\}$ coincide in a finite number of steps, and their convergence is equivalent.
In the rest of this section,   we   demonstrate   in detail   the aforementioned ideas.

\subsection{Auxiliary Sequences}

Denote by $\tau_{i,m}\deq\textrm{inf}  \{k: \sigma_{i,k}=m \}$
   the  smallest time  when   the truncation number of agent  $i$ has reached    $m$,  and by
          $\tau_m\deq \min\limits_{i \in \mathcal{V}}\tau_{i,m} $
 the smallest time  when at least one  agent  has its truncation number reached  $  m$.
 Set     $\tilde{\tau}_{j,m}\deq \tau_{j,m} \wedge \tau_{m+1} $, where
 $a\wedge b =\min\{a,b\}. $

For any $i \in \mathcal{V}$,   we  construct two auxiliary sequences $ \{\tilde{x}_{i,k}\}_{k \geq 0}$   and $ \{\tilde{\varepsilon}_{i,k+1} \}_{k \geq 0}$.
Note that for  the  considered    $\omega$,  for  any  $k \geq 0$ there exists a unique   integer  $m \geq 0$ such that $ k\in [\tau_m , \tau_{m+1}).$  
We  then define
 \begin{align}
    & \tilde{x}_{i,k}\triangleq  x^{*}, ~  \tilde{\varepsilon}_{i,k+1}\triangleq  -f_i(x^{*})\quad {\rm if~}    k  \in[ \tau_m,   \tilde{\tau}_{i,m}),  \label{aux0} \\
     &  \tilde{x}_{i,k}\triangleq  x_{i,k}, ~  \tilde{\varepsilon}_{i,k+1}\triangleq  \varepsilon_{i,k+1}  \quad {\rm if~}   k \in [ \tilde{\tau}_{i,m} , \tau_{m+1}) .     \label{aux1}
\end{align} 
  \begin{lem}\label{lemau2}
The sequences   $ \{\tilde{x}_{i,k}\}, \{ \tilde{\varepsilon}_{i,k+1}\} $ defined by \eqref{aux0} \eqref{aux1}  satisfy the following recursions
\begin{align}
       \hat{x}_{i,k+1} & \deq\sum_{j  \in N_i(k)}  \omega_{ij} (k)    \tilde{x}_{j,k} +
    \gamma_{ k } ( f_i(  \tilde{x}_{i,k})+\tilde{\varepsilon}_{i,k+1}),  \label{algorithm01} \\
     \tilde{x}_{i,k+1}  &  =\hat{x}_{i,k+1} I_{  \left[  \parallel \hat{x}_{j,k+1}  \parallel \leq M_{\sigma_k}~\forall  j \in \mathcal{V} \right] } \nonumber
 \\& +x^* I_{ \left[    \exists   j \in \mathcal{V} ~  \parallel  \hat{x}_{j,k+1} \parallel > M_{\sigma_k}   \right]}, \label{algorithm02}  \\
   \sigma_{k+1} &=\sigma_k+ I_{ \left[ \exists   j \in \mathcal{V}  ~\parallel  \hat{x}_{j,k+1} \parallel > M_{\sigma_k} \right]}, \quad   \sigma_0=0.    \label{algorithm03}
\end{align}
\end{lem}

The proof is given  in Appendix \ref{PLA}.

Before clarifying  the property of the noise sequence   $ \{\tilde{\varepsilon}_{i,k+1} \}_{k \geq 0}$,
we  need   the following lemma.
 \begin{lem}\label{lem001}
 Assume A4 holds. Then
 \begin{flalign}
&~~\textrm{ i)}  ~~~  \sigma_{j,k+B d_{i,j}} \geq   \sigma_{i,k}~~~\forall j \in \mathcal{V}~~ \forall k \geq 0,   &  \label{gap000}
\end{flalign}
  where $d_{i,j}$  is the length   of  the shortest  directed path from  $i$ to $j $ in   $\mathcal{G}_{\infty}$,
     and $B $ is   the positive integer given  in A4 d),
      \begin{flalign}
&~~~\textrm{ ii)}  ~~~ \tilde{\tau}_{j,m} \leq \tau_m+BD \quad \forall j \in \mathcal{V} ~\textrm{ for } m \geq 1,  & \label{gap11}
\end{flalign}
 where  $D\deq\max\limits_{i,j \in \mathcal{V}}d_{i,j}$.
\end{lem}
 
\begin{IEEEproof}
i)
Since $\mathcal{G}_{\infty}$  is strongly connected by A4  c), for any $j \in \mathcal{V}$ there exists a sequence of nodes  $i_1, i_2, \cdots, i_{d_{i,j}-1}$  such that  $(i, i_1) \in \mathcal{E}_{\infty},(i_1, i_2) \in \mathcal{E}_{\infty},
 \cdots,  (i_{d_{i,j}-1},j)   \in \mathcal{E}_{\infty}$.

Noticing  that   $ ( i,i_1) \in  \mathcal{E}_{\infty},$   by A4  d) we have
\begin{equation}
(i,i_1) \in \mathcal{E}(k)\cup \mathcal{E}(k+1)  \cup  \cdots \cup \mathcal{E}(k+B-1). \nonumber
\end{equation}
  Therefore,  there exists a positive integer $k' \in [k ,k+B-1]$ such that
$(i,i_1) \in \mathcal{E}(k') .$
So,  $i \in N_{i_1}(k')$,  and hence by   \eqref{comp1} and \eqref{step2}  we derive
 \begin{equation}\label{ineq}
  \begin{split}
    \sigma_{i_1,k+B}  & \geq \sigma_{i_1,k'+1} \geq   \hat{\sigma}_{i_1,k'}   \geq   \sigma_{i,k'}\geq   \sigma_{i,k } . \nonumber
    \end{split}
 \end{equation}
Similarly, we have
$\sigma_{i_2,k+2B} \geq  \sigma_{i_1,k+B} \geq  \sigma_{i,k } .$
Continuing this  procedure, we finally reach the   inequality \eqref{gap000}.

  ii)
Let $\tau_m=k_1 $  for some  $m\geq1$.  Then  there is an $i$ such that $\tau_{i,m}=k_1.$
By \eqref{gap000} we have
$
  \sigma_{j,k_1+B d_{i,j}} \geq  \sigma_{i,k_1 }=m  ~~\forall j \in \mathcal{V}. $

For the case where $\sigma_{j,k_1+ Bd_{i,j}} = m  ~\forall j \in \mathcal{V},$ we have
$\tau_{j,m}\leq k_1+Bd_{i,j}  ~\forall j \in \mathcal{V}.$ By noticing $\tau_m=k_1,$  from here
 by the definition of $\tilde{\tau}_{j,m}$  we obtain  \eqref{gap11}:
 \begin{align*}
\tilde{\tau}_{j,m}\leq \tau_{j,m}\leq \tau_m+Bd_{i,j}\leq \tau_m+BD ~~\forall j \in \mathcal{V}.
\end{align*}
For the case where $\sigma_{j,k_1+ Bd_{i,j}} > m$ for some $j \in \mathcal{V},$  we     have
$\tau_{m+1}\leq k_1+Bd_{i,j} \textrm{ for some }j \in \mathcal{V}, $ and hence  $\tau_{m+1}\leq  \tau_m+BD $.
Again, by noticing $\tau_m=k_1 $  we obtain \eqref{gap11}:
 \begin{align*}
\tilde{\tau}_{j,m}\leq \tau_{m+1}\leq   \tau_m+BD ~~\forall j \in \mathcal{V}.
\end{align*}
\end{IEEEproof}

\begin{prp} \label{cor}  If $  \sigma_k \xlongrightarrow [k \rightarrow \infty ]{} \infty$,  then for any  $i \in \mathcal{V},$
$ \lim\limits_{k \rightarrow \infty } \sigma_{i,k} = \infty$.
This is because there exists  an $i_0 \in \mathcal{V}$ such that $ \sigma_{i_0,k} \xlongrightarrow [k \rightarrow \infty ]{} \infty$.
Then from \eqref{gap000} it follows that for any $j \in \mathcal{V},$  $\sigma_{j,k} \xlongrightarrow [k \rightarrow \infty ]{} \infty $.
\end{prp}

\begin{prp}\label{r3} If $\{ \sigma_k\}$ is bounded, then $\{ \tilde{x}_{i,k}\}$   and  $\{x_{i,k}\}$, $\{ \tilde{\varepsilon}_{i,k}\}$   and  $\{\varepsilon_{i,k}\}$  coincide
in a finite number of steps.   \end{prp}
\begin{IEEEproof}Since    $\sigma_{k}  $     is defined as
 the largest truncation number among all agents at time $k$,  from $\lim\limits_{k \rightarrow \infty} \sigma_{k} =  \sigma$   it follows that   for any $i \in \mathcal{V}$ and $k\geq 0,$ $ \sigma_{i, k} \leq \sigma  .$
Then by the definition of $\tau_{i,m}$ we have that for any $i \in \mathcal{V},$
 $\tau_{i, \sigma+1} = \textrm{inf}  \{k: \sigma_{i,k}= \sigma+1 \} = \infty  ,$ and hence 
  \begin{equation}\label{infinity}
    \tau_{\sigma+1}=\infty\textrm{  when } \lim\limits_{k \rightarrow \infty} \sigma_{k} =  \sigma .
\end{equation}
Then the  result is derived  by  \eqref{aux1}. \end{IEEEproof}
 \begin{lem}\label{lemma1}
Assume  A5  b) holds at  the  sample path  $\omega$ under consideration for all agents.   Then for this $\omega$
\begin{small}
  \begin{equation}\label{auxnoise}
  \begin{split}
   & \lim\limits_{T \rightarrow 0} \limsup\limits_{k \rightarrow \infty} \frac{1}{T}
\left \|  \sum_{s=n_k}^{m(n_k,t_k)\wedge ( \tau_{\sigma_{n_k}+1}-1 )} \gamma_s  \tilde{\varepsilon}_{ s+1}
    I_{[ \parallel \tilde{X}_{  s}\parallel \leq K ]} \right \|
    \\&~~~~~~~    =0 ~~\forall t_k \in[0,T] \textrm{   for sufficiently large  } K  >0
    \end{split}
  \end{equation}
  \end{small}
  along indices $\{ n_k\}$
  whenever  $\{ \widetilde{X}_{  n_k}\}$   converges  at    $\omega$.   \end{lem}

\begin{IEEEproof}  The proof is shown in Appendix \ref{PL1}.
\end{IEEEproof}

\subsection{Local  Properties  Along   Convergent Subsequences}

  Set  $\Psi(k,s)\triangleq  [ D_{\bot} (W(k) \otimes \mathbf{I}_l)]  [ D_{\bot} (W(k-1) \otimes \mathbf{I}_l)] \cdots$
  $(W(s) \otimes \mathbf{I}_l)]~~\forall k\geq s,  \quad  \Psi(k-1,k)\triangleq   \mathbf{ I}_{Nl}   .$

  Beacuse  for any $k\geq 1,$ the   matrix $  W(k )  $ is doubly stochastic, then
  by using the rule of Kronecker product \begin{equation}\label{kr}
(A \otimes B)(C \otimes D)= (AC ) \otimes ( BD),
\end{equation}  we  conclude that for any $ k \geq s-1$,
  \begin{align}
      & \Psi(k,s)=   \left (\Phi(k,s)  - \frac{1}{N} \mathbf{1} \mathbf{1}^T \right) \otimes \mathbf{I}_l ,\label{power23}\\
&  \Psi(k,s) D_{\bot} =   \left( \Phi(k,s)  - \frac{1}{N} \mathbf{1} \mathbf{1}^T \right)  \otimes \mathbf{ I}_l  . \label{power24}
\end{align}

The  following  lemma    measures the closeness of
the auxiliary  sequence  $\{ \widetilde{X}_{ k}\}$ along its   convergent subsequence  $\{  \widetilde{X}_{n_k}\}$.
   \begin{lem}\label{lemma2}
    Assume A1, A3, A4 hold and  that   A5  b) holds at the sample path $\omega$  under consideration.
 Let   $\{ \widetilde{X}_{n_k}\}$  be a convergent subsequence of $\{ \widetilde{X}_k\}:
 \widetilde{X}_{n_k} \xlongrightarrow [k \rightarrow \infty] {} \bar{X}$
 at the considered  $\omega$.  Then for this $\omega$, there is a  $T>0$
 such that for sufficiently large   $k$ and any  $T_k \in [0,T]$,
   \begin{equation}\label{compactform}
 \widetilde{X}_{m+1}= \left( W(m) \otimes \mathbf{I}_l\right)
   \widetilde{X}_m+ \gamma_{m } \left(F(\widetilde{X}_m)+\tilde{\varepsilon}_{m+1}\right)
  \end{equation}
 for any $m=n_k, \cdots, m(n_k, T_k)$,    and
         \begin{align}
   &   \parallel  \widetilde{X}_{m+1} -\widetilde{X}_{n_k}  \parallel \leq c_1 T_k +M_0',  \label{res21}\\
   &   \parallel  \bar{x}_{m+1}-\bar{x}_{n_k}  \parallel \leq c_2 T_k~  ~ \forall m: n_k \leq m \leq m(n_k,T_k) , \label{res22}
\end{align}
 where     $ \bar{x}_k\triangleq \frac{1}{N} ( \mathbf{1}^T \otimes \mathbf{I}_l)  \widetilde{X}_k =\frac{1}{N}\sum_{i=1} ^N \tilde{x}_{i,k} $
 and   $c_0,~c_1 ,~M_0' $ are   positive  constants  which  may   depend  on $\omega$.
  \end{lem}

The proof   is given in  Appendix \ref{PL2}.

Multiplying both sides of \eqref{compactform} with $ \frac{1}{N} ( \mathbf{1}^T \otimes \mathbf{I}_l) $
 from left, by $ \mathbf{1}^TW(m) =\mathbf{1}^T$  and   \eqref{kr}
it follows that   \begin{equation}\label{centrelized}
  \begin{split}
   \bar{x}_{m+1}= \bar{x}_m &+ \gamma_m f(\bar{x}_m)+
  (\mathbf{1}^T \otimes \mathbf{I}_l)\gamma_m  \tilde{\varepsilon}_{m+1}/ N \\&+
      \gamma_m  \sum_{i=1}^N  \left(  f_i(\tilde{x}_{i,m})- f_i(\bar{x}_m) \right)/N.    \end{split}
\end{equation}
 By setting  $e_{i,m+1}\triangleq \left(  f_i(\tilde{x}_{i,m})- f_i(\bar{x}_m) \right)/N $,
  $e_{m+1}\triangleq \sum_{i=1}^N e_{i,m+1}$, and
$  \zeta_{m+1}\triangleq  (\mathbf{1}^T \otimes \mathbf{I}_l)   \tilde{\varepsilon}_{m+1}/N+  e_{m+1},$
 we  can rewrite \eqref{centrelized} in the centralized form   as follows:
\begin{equation}\label{centrelized2}
  \begin{split}
  \bar{x}_{m+1}=\bar{x}_m +  \gamma_m \left( f(\bar{x}_m)+  \zeta_{m+1} \right).
    \end{split}
\end{equation}

The following lemma gives  the property of   the noise sequence $\{  \zeta_{k+1}\}$.
For the  proof  we refer to  Appendix \ref{PLN}.

  \begin{lem}\label{lemnoise}
  Assume  that all  conditions used  in  Lemma  \ref{lemma2} are satisfied.
 Let   $\{ \widetilde{X}_{n_k}\}$  be a convergent subsequence with limit $\bar{X}$ at the considered  $\omega$.  Then for this $\omega$
\begin{equation}\label{err2}
      \lim_{T \rightarrow 0} \limsup_{k \rightarrow \infty} \frac{1}{T}
  \left \|  \sum_{s=n_k}^{m(n_k,T_k) } \gamma_s   \zeta_{s+1}  \right\|=0 \quad \forall T_k \in[0,T] .
    \end{equation}
\end{lem}

   The following lemma   gives  the  crossing behavior of $v(\cdot)$  at
 the trajectory $\bar{x}_k$   with respect to a non-empty interval that  has  no intersection with
 $v(J)$.

   \begin{lem}\label{lemma3}
  Assume A1-A4 hold and that A5 b) holds   for the sample path $\omega$ under consideration.
   Then any nonempty interval $[\delta_1, \delta_2 ]$ with    $d([\delta_1, \delta_2 ],  v(J))>0 $
  cannot be crossed by    $\{v(\bar{x}_{n_k}), \cdots, v(\bar{x}_{m_k})\}$  infinitely many times  with
     $ \{   \| \widetilde{X}_{n_k}  \|  \}$  bounded,
   where by `` crossing  $[\delta_1, \delta_2]$ by  $\{v(\bar{x}_{n_k}), \cdots,  v(\bar{x}_{m_k})\}$ "
   it is meant that   $v(\bar{x}_{n_k}) \leq \delta_1,  v(\bar{x}_{m_k}) \geq \delta_2 ,$
and $    \delta_1  < v(\bar{x}_s) < \delta_2 ~~\forall s: n_k< s<m_k$.
  \end{lem}
 \begin{IEEEproof}
  Assume the converse: for some nonempty interval $[\delta_1, \delta_2 ]$ with
    $d([\delta_1, \delta_2 ],  v(J))>0 $, there exist  infinitely many crossings
  $\left \{ v(\bar{x}_{n_k}), \cdots, v(\bar{x}_{m_k}) \right \}$  with $ \{     \|\widetilde{X}_{n_k}\| \} $  bounded.

  By the boundedness of  $ \{  \|\widetilde{X}_{n_k}\|  \}$, we  can extract a  convergent subsequence
   still denoted by  $\{ \widetilde{X}_{n_k}\}$ with
    $\lim\limits_{k \rightarrow \infty} \widetilde{X}_{n_k} =  \bar{X}$. So,
     $\lim\limits_{k \rightarrow \infty}  \bar{x}_{n_k}= \bar{x}= \frac{\mathbf{1}^T \otimes \mathbf{I}_l}{N} \bar{X}$.  Setting  $T_k=\gamma_{n_k}$ in \eqref{res22}, we derive
    $ \parallel  \bar{x}_{n_k+1}-\bar{x}_{n_k}  \parallel \leq c_2 \gamma_{n_k}
    \xlongrightarrow [k \rightarrow \infty]{} 0   .$
By the definition of crossings,     we obtain    $v(\bar{x}_{n_k}) \leq \delta_1 < v(\bar{x}_{n_k+1})$, and hence
   \begin{equation}\label{gap1}
  v(\bar{x}_{n_k})  \xlongrightarrow [k \rightarrow \infty]{} \delta_1 =v(\bar{x}),
  \quad  d(\bar{x}, J) \triangleq  \vartheta>0.
   \end{equation}
 Then by   \eqref{res22} it is seen  that
   \begin{equation}\label{gap2}
    d(\bar{x}_s, J)  > \frac{\vartheta}{2}  \quad \forall s : n_k \leq s \leq m(n_k,t)+1
   \end{equation}
for sufficiently small  $t>0$  and large $k$.
  From  \eqref{centrelized2} we obtain
\begin{equation*}
  \begin{array}{lll}
   &    v(\bar{x}_{m(n_k,t)+1})= v(\bar{x}_{n_k} + \sum_{s=n_k}^{m(n_k,t)} \gamma_s(   f(\bar{x}_s)+  \zeta_{s+1})) \\
   &  = v(\bar{x}_{n_k})+ v_x(\xi_k)^T  \sum_{s=n_k}^{m(n_k,t)} \gamma_s( f(\bar{x}_s)+  \zeta_{s+1}),
      \end{array}
\end{equation*}
where  $\xi_k$ is  in-between $\bar{x}_{n_k}$ and $\bar{x}_{m(n_k,t)+1 } $. Then
\begin{small}
\begin{equation}\label{vfunction2}
  \begin{array}{lll}
   & v(\bar{x}_{m(n_k,t)+1})  - v(\bar{x}_{n_k})    =  \sum_{s=n_k}^{m(n_k,t)}\gamma_s v_x( \bar{x}_s)^T  f(\bar{x}_s)
  \\
   &  +  \sum_{s=n_k}^{m(n_k,t)} \gamma_s \left (v_x(\xi_k)-   v_x( \bar{x}_s) \right)^T   f(\bar{x}_s)    \\& +v_x(\xi_k)^T  \sum_{s=n_k}^{m(n_k,t)} \gamma_s  \zeta_{s+1}.            \end{array}
\end{equation}
\end{small}

By   \eqref{4.1} and   \eqref{gap2},     there exists  a constant $\alpha_1>0$ such that
$$ v_x( \bar{x}_s)^T  f(\bar{x}_s) \leq  - \alpha_1 \quad \forall s : n_k \leq s \leq m(n_k,t) $$
 for  sufficiently small  $t>0$  and large $k$, and hence
 \begin{equation}\label{estm41}
  \begin{array}{lll}
   \sum_{s=n_k}^{m(n_k,t)}  \gamma_s  v_x( \bar{x}_s)^T  f(\bar{x}_s) \leq  - \alpha_1 t.
   \end{array}
\end{equation}

Since    $\{\bar{x}_{s } : n_k \leq s \leq  m(n_k,t) \}$ are bounded,   by continuity of $f(\cdot)$  there exists a constant  $c_6>0$ such that
 \begin{equation}\label{sum41}
  \begin{array}{lll}
  \sum_{s=n_k}^{m(n_k,t)} \gamma_s \parallel f(\bar{x}_s) \parallel \leq  c_6t .
   \end{array}
\end{equation}
By invoking that  $\xi_k$ is  in-between $\bar{x}_{n_k}$ and $\bar{x}_{m(n_k,t)+1 } $,
 by continuity of  $v_x(\cdot)$ and   \eqref{res22},  we know that
 \begin{equation}\label{D01}
    v_x(\xi_k)
-  v_x( \bar{x}_s) =\delta(t)   \quad \forall s: n_k \leq s\leq m(n_k,t),
 \end{equation}
  where $\delta(t) \rightarrow 0\textrm{ as }t \rightarrow 0$.
Then by  \eqref{sum41}   we derive
  \begin{equation}\label{sum42}
  \begin{array}{lll}
      & \sum_{s=n_k}^{m(n_k,t)} \gamma_s (v_x(\xi_k)-   v_x( \bar{x}_s))^T  f(\bar{x}_s)
     \\& \leq    \delta(t) \sum\limits_{s=n_k}^{m(n_k,t)} \gamma_s \parallel  f(\bar{x}_s) \parallel \leq \delta(t) t .
     \end{array}
  \end{equation}

 Since  $ \bar{x}_{n_k} \xlongrightarrow [k \rightarrow \infty]{}  \bar{x}$,    by continuity of  $v_x(\cdot)$  and   \eqref{res22} it follows  that for  sufficiently small  $t>0$  and large $k$
  $$  v(\bar{x}_{s})-v_x( \bar{x} )= o(1)+\delta(t)   \quad \forall s: n_k \leq s \leq m(n_k,t) ,$$
    where $o(1)\rightarrow 0\textrm{ as } k \rightarrow \infty.$ Then by \eqref{D01}  we derive
 $$ v_x(\xi_k)-v_x( \bar{x} )= o(1)+\delta(t) \quad \forall s: n_k \leq s \leq m(n_k,t).$$
Consequently, for  sufficiently small  $t>0$  and large $k$
    \begin{equation}\label{sum43}
    \begin{array}{lll}
    & v_x(\xi_k)^T   \sum_{s=n_k}^{m(n_k,t)} \gamma_s   \zeta_{s+1} \\&=\left( (v_x(\xi_k)-v_x( \bar{x} )) + v_x( \bar{x} ) \right)^T  \sum_{s=n_k}^{m(n_k,t)} \gamma_s  \zeta_{s+1}  \\
& \leq  (o(1)+\delta(t)+\|v_x( \bar{x} )\|) \| \sum_{s=n_k}^{m(n_k,t)} \gamma_s       \zeta_{s+1} \|.
    \end{array}
  \end{equation}
Substituting   \eqref{estm41},  \eqref{sum42}, and  \eqref{sum43} into   \eqref{vfunction2}, we obtain
  \begin{equation*}
    \begin{array}{lll}
  &   v(\bar{x}_{m(n_k,t)+1}) - v(\bar{x}_{n_k}) \leq - \alpha_1 t  +  \delta(t) t
 \\& + (o(1)+\delta(t)+ \| v_x( \bar{x} ) \|) \| \sum\limits_{s=n_k}^{m(n_k,t)} \gamma_s       \zeta_{s+1} \|.
    \end{array}
  \end{equation*}
Then  by   \eqref{gap1}  it follows that
\begin{equation*}
   \begin{array}{lll}
  &  \limsup\limits_{k \rightarrow \infty} v(\bar{x}_{m(n_k,t)+1}) \leq \delta_1 - \alpha_1 t+    \delta(t) t \\& +
 \left(  \delta(t)  +\| v_x( \bar{x} ) \| \right)     \limsup\limits_{k \rightarrow \infty} \|  \sum_{s=n_k}^{m(n_k,t)}
  \gamma_s   \zeta_{s+1}\|.
 \end{array}
\end{equation*}
Hence from \eqref{err2},     for sufficiently small  $t$ we have that
\begin{equation}\label{vfunction4}
  \begin{split}
   &  \limsup_{k \rightarrow \infty}  v(\bar{x}_{m(n_k,t)+1})  \leq \delta_1- \frac{\alpha_1}{2} t .      \end{split}
\end{equation}

 However,     by continuity of  $v_x(\cdot)$ and    \eqref{res22} we know that
$$ \lim_{t \rightarrow 0}  \max_{n_k \leq m \leq m(n_k,t)}\parallel v(\bar{x}_{m+1})-
  v(\bar{x}_{n_k}) \parallel=0,
$$
 which implies that   $m(n_k,t) +1 < m_k$ for sufficiently small $t$. Therefore,   $ v(\bar{x}_{m(n_k,t)+1})  \in  (\delta_1, \delta_2 )$, which  contradicts   with \eqref{vfunction4}.
 Consequently, the converse assumption   is not true. The proof is completed.
\end{IEEEproof}

\subsection{Finiteness of Number of Truncations}
     \begin{lem}\label{lemma4}
Let $ \{ x_{i,k}\} $ be produced by \eqref{comp1}-\eqref{step2}  with an arbitrary  initial value  $   x_{i,0} $.
Suppose  A1, A3, A4, and assume A5 b) holds for the sample path $\omega$ under consideration.

i) If $\lim\limits_{k \rightarrow \infty } \sigma_k =  \infty, $  then there exists an integer sequence $\{n_k\} $
 such that     $\bar{x}_{n_k}= x^{*} $. Further,  $\{\bar{x}_{n_k}\}$   starting from $ x^{*}$  crosses
   the sphere with   $\| x \| =c_0 $  infinitely many times, where  $\{\bar{x}_{k}\}$  is defined in Lemma 4.6 and $c_0$ is given in A2 c).

ii) If, in addition, A2 also holds, then
 there exists a positive  integer  $\sigma$ possibly depending on $\omega$ such   that
 \begin{equation}\label{lem4}
 \lim_{k \rightarrow \infty} \sigma_{ k}=\sigma  .
\end{equation}
\end{lem}
The proof of  the lemma is given in Appendix \ref{PL4}.
Lemma  \ref{lemma4} says    that  the largest   truncation number  among all agents  converges, while
the   following lemma   indicates  that the  truncation numbers at  all  agents    converge  to the  same limit.
 \begin{lem}\label{lem001}
 Assume    all  conditions required by  Lemma  \ref{lemma4}  are satisfied.
  Then  there exists a positive  integer  $\sigma$ such   that
\begin{equation}\label{trb}
 \lim_{k \rightarrow \infty} \sigma_{i, k}=\sigma ~~ \forall i\in \mathcal{V}.
  \end{equation}
\end{lem}

\begin{IEEEproof}   Since  all  conditions required by Lemma  \ref{lemma4}  hold,   \eqref{lem4} holds for some  positive  integer  $\sigma$.
 Thus,
 \begin{equation}\label{upper}
 \sigma_{i, k} \leq \sigma ~~ \forall k \geq 0 ~~ \forall i \in \mathcal{V}.
\end{equation}
By \eqref{lem4}  and \eqref{infinity},  we have  $\tau_{\sigma+1}=\infty$, and hence  $\tilde{\tau}_{i,\sigma}  =\tau_{i,\sigma} \leq BD+ \tau_{\sigma} ~ \forall i \in \mathcal{V} $ by \eqref{gap11}.
This means  that the smallest time when  the truncation number
 of agent $i$ reaches $\sigma$ is not larger than $ BD+ \tau_{\sigma}$.
So, the truncation number  of  agent $i$ after time $ BD+ \tau_{\sigma}$ is not smaller than  $\sigma$,
i.e.,  $  \sigma_{i, k} \geq \sigma ~ \forall k \geq BD+ \tau_{\sigma} ~  \forall i \in \mathcal{V}.$
It incorporating  with \eqref{upper}  yields \begin{equation}\label{consensus}
 \sigma_{i, k}= \sigma ~~ \forall k \geq BD+ \tau_{\sigma} ~~ \forall i \in \mathcal{V}.
\end{equation}
 Consequently, \eqref{trb} holds.
\end{IEEEproof}

\subsection{  Proof  of Theorem  \ref{thm1}}
{\bf Proof.} i) By   \eqref{lem4} and \eqref{consensus},
 there is a positive integer $\sigma$  possibly  depending  on $\omega$ such  that
 \begin{equation}\label{hat}
\hat{\sigma}_{i, k}=\sigma_{i, k}= \sigma ~~ \forall k \geq k_0 \deq BD+ \tau_{\sigma} ~~ \forall i \in \mathcal{V}.
\end{equation}
Then by \eqref{comp2}, we have the following for all $  k \geq k_0 $:
 \begin{equation}\label{prime}
 x'_{i,k+1}=  \sum_{ j \in  N_{i }(k)}  \omega_{ij}(k)   x_{j,k}   + \gamma_{ k } O_{i,k+1}
~~   \forall i \in \mathcal{V}.
\end{equation}
 By \eqref{consensus} and \eqref{hat},    we obtain $\sigma_{i,k+1}=\hat{\sigma}_{i,k} = \sigma  ~ \forall k \geq k_0~\forall i \in \mathcal{V}$.
  Thus,    for any $ k \geq k_0$ and any $i \in \mathcal{V}$,  we have that $ \|x_{i,k+1}' \| \leq M_{\sigma}  $ by  \eqref{step2},   and
  $ x_{i,k+1}=x_{i,k+1}' $  by \eqref{step1}. Then    for any  $i \in \mathcal{V}$,     $\{x_{i,k}\}$ is bounded and \eqref{notr}  follows from \eqref{prime}.

ii) By multiplying both sides of  \eqref{centalform} with $ D_{\bot}$   from left, we have that
  \begin{equation}
   \begin{split}
  & X_{\bot, k+1}=  D_{\bot} (W(k) \otimes I_l)  X_{  k }+ \gamma_{k }  D_{\bot}     (F( X_k)+\varepsilon_{k+1})  ,  \nonumber
  \end{split}
   \end{equation}
 and inductively
   \begin{equation}\label{dis021}
   \begin{array}{lll}
  & X_{\bot, k+1}= \Psi(k, k_0) X_{k_0}+ \\
  &   \sum_{m=k_0} ^k
  \gamma_{m } \Psi(k-1, m) D_{\bot} (F( X_m)+ \varepsilon_{m+1}) \quad \forall k \geq  k_0. \nonumber
  \end{array}
   \end{equation}
  Then by   \eqref{power23} and \eqref{power24},  we  derive
    \begin{equation}\label{disagree}
   \begin{array}{lll}
  & X_{\bot, k+1}   =  [(\Phi(k, k_0) - \frac{1}{N} \mathbf{1} \mathbf{1}^T) \otimes  \mathbf{I}_l]    X_{k_0}\\&+ \sum_{m=k_0} ^k
  \gamma_{m }   [(\Phi(k-1, m) - \frac{1}{N} \mathbf{1} \mathbf{1}^T) \otimes \mathbf{I}_l]  F(X_m) \\
  & +\sum_{m=k_0} ^k \gamma_{m }    [( \Phi(k-1, m) - \frac{1}{N} \mathbf{1} \mathbf{1}^T) \otimes \mathbf{I}_l]    \varepsilon_{m+1}.  \notag
  \end{array}
   \end{equation}
Therefore,  from  \eqref{graphmatrix}  by continuity of $F(\cdot )$ and
   the boundedness of $\{ X_s\}$, we  conclude  that there exist positive  constants $c_1 , c_2  ,c_3  $
    possibly depending  on $\omega$   such that
  \begin{equation}\label{er2}
   \begin{array}{lll}
  & \parallel X_{\bot, k+1} \parallel \leq  c_1  \rho^{k+1- k_0} + c_2  \sum_{m= k_0} ^k \gamma_m \rho^{k -m}    \\
  & +  c_3\sum_{m=k_0} ^{k} \gamma_m  \rho ^{k-m}   \parallel \varepsilon _{m+1}\parallel
  \quad \forall k \geq  k_0 .
  \end{array}
   \end{equation}

   Note  that for any given  $\epsilon>0$, there exists a positive integer $k_1 $ such that
$ \gamma_k \leq \epsilon~~ \forall k \geq k_1$. We then  have
 \begin{equation}
\begin{array}{lll}
&  \sum_{m=0}^k \gamma_m \rho^{k-m} = \sum_{m=0}^{k_1}  \gamma_m \rho^{k-m}
+ \sum_{m=k_1+1}^{k}  \gamma_m \rho^{k-m}
 \\&  \leq \rho^{k-k_1}   \sum_{m=0}^{k_1 }  \gamma_m +  \epsilon \frac{1}{1-\rho} \xlongrightarrow [k \rightarrow \infty \atop \epsilon \rightarrow 0  ]{}  0  \nonumber.
\end{array}
\end{equation}
Therefore, the second  term  at the right-hand side of \eqref{er2} tends  to zero as $k \rightarrow \infty$.
Similarly,  the last term  at the right-hand side of \eqref{er2} also tends to zero  since
 $ \lim\limits_{k \rightarrow \infty} \gamma_k  \varepsilon_{k+1}=0  .$
 Therefore, by $0< \rho<1$, from \eqref{er2} we conclude that
    $$  X_{\bot, k} \xlongrightarrow [k \rightarrow \infty]{} \mathbf{0}.$$

By   i) and Proposition    \ref{r3},    we see  that  $\{\tilde{x}_{i,k}~\forall i\in \mathcal{V}\} $
   are bounded for this  $\omega$, and hence $\{\bar{x}_k\}$ is bounded.
   The rest of the  proof  is similar to that    given in \cite{Chen_2002}.

We first show the convergence of $ v( \bar{x}_k).$ Since
 $$ v_1 \triangleq \liminf_{k \rightarrow \infty } v( \bar{x}_k) \leq  \limsup_{k \rightarrow \infty } v( \bar{x}_k)
  \triangleq v_2 ,$$
we want to prove $v_1=v_2$.
Assume the converse: $v_1< v_2$. Since   $v(J)$  is nowhere dense, there exists a nonnegative interval   $[\delta_1, \delta_2 ]
    \in ( v_1, v_2)$ such that    $d([\delta_1, \delta_2 ],  v(J))>0 $.
Then  $v(\bar{x}_k)$ crosses the interval   $[\delta_1, \delta_2 ]$ infinitely  many times.
This contradicts  Lemma \ref{lemma3}.
Therefore, $v_1=v_2$, which implies the convergence of  $v( \bar{x}_k)$.

We  then prove $ d(  \bar{x}_k ,J )  \xlongrightarrow [k \rightarrow \infty]{} 0$.
Assume the converse. Then by the boundedness of  $ \{\bar{x}_k\}$,   there   exists a convergent subsequence
$  \bar{x}_{n_k} \xlongrightarrow [k \rightarrow \infty]{} \bar{x} $ with
 $ d(\bar{x}, J) \triangleq  \vartheta>0$. From  \eqref{res22} it follows  that for
sufficiently   small    $t>0$  and large  $k$
   $$  d(\bar{x}_s, J)  > \frac{\vartheta}{2} \quad \forall s : n_k \leq s \leq m(n_k,t), $$
and hence by  \eqref{4.1}   there exists a  constant  $b>0$ such that
   $$  v_x(\bar{x}_s)^T f(\bar{x}_s) < -b   \quad \forall s: n_k \leq s \leq m(n_k,t).$$
   Thus, similar to the proof  for obtaining  \eqref{vfunction4} it is seen  that for sufficiently small $t>0 $
   \begin{equation} \label{contradict1}
   \limsup_{k \rightarrow \infty}  v(\bar{x}_{m(n_k,t)+1})  \leq  v(\bar{x})- \frac{b}{2} t  \nonumber.
\end{equation}
This contradicts with the convergence of  $v( \bar{x}_k)$.
Therefore,  $ d(  \bar{x}_k ,J )  \xlongrightarrow [k \rightarrow \infty]{} 0$, and hence
  $ d( x_k ,J )  \xlongrightarrow [k \rightarrow \infty]{} 0$.

iii)
Assume the converse: i.e., $J^{*}$ is disconnected. Then there exist closed sets $J_1^{*}$ and
$J_2^{*}$ such that $J^{*}= J_1^{*} \cup J_2^{*}$ and $d(J_1^{*}, J_2^{*})>0$.
Define $\rho=\frac{1}{3}d(J_1^{*}, J_2^{*}).$
 By   $ d( \bar{x}_k ,J^{*} )  \xlongrightarrow [k \rightarrow \infty]{} 0$  there exists $k_0$ such that
\begin{equation}
\label{E1}
 \bar{x}_k \in B(J_1^{*}, \rho) \cup   B(J_2^{*}, \rho) ~~ \forall k\geq k_0,
\end{equation}
where $B(A, \rho)$ denotes the $\rho$-neighborhood   of the set $A.$

Define
\begin{align}
& n_0=\inf \{ k >k_0, d( \bar{x}_k ,J_1^{*} )<\rho \},  \nonumber\\
& m_p=\inf \{ k >n_p, d( \bar{x}_k ,J_2^{*} )<\rho \}, \nonumber\\
& n_p=\inf \{ k >m_p, d( \bar{x}_k ,J_1^{*} )<\rho \},   ~~ p\geq 0 .\nonumber
\end{align}
 Then by \eqref{E1},  we have that   $  \bar{x}_{n_p} \in B(J_1^{*}, \rho)$ and  $\bar{x}_{n_p+1} \in B(J_2^{*}, \rho)  $
   for any $ p \geq 0$. Then by   $d(J_1^{*}, J_2^{*})=3\rho$,  we derive
\begin{equation} \label{E2}
 \|  \bar{x}_{n_p} - \bar{x}_{n_p+1} \| >\rho .
\end{equation}

Since $\{\bar{x}_{n_k} \}$ is bounded, we can extract a convergent subsequence,
denoted  still  by  $\{\bar{x}_{n_k} \}$.   By setting $T_k = \gamma_{n_k }$ in  \eqref{res22}  we derive
 $   \|  \bar{x}_{n_k+1} - \bar{x}_{n_k } \| \leq c_2 \gamma_{n_k }   \xlongrightarrow [k \rightarrow \infty]{} 0, $
which contradicts with \eqref{E2}.   So,
 the converse assumption is not true. The  proof is completed.
\hfill $\blacksquare$

\section{Convergence for Application Problems}
In this section, we establish   convergence results  for the two problems stated  in Section II,
and present the corresponding numerical  simulations.

\subsection{Distributed PCA}
 We now apply DSAAWET to  the problem of distributed principal component analysis,
and establish  its convergence.
\begin{cor}\label{thm-PCA}Let  $\{x_{i,k}\}$ be produced by \eqref{comp1}--\eqref{step2}, where $O_{i,k+1}$ is defined by \eqref{pca}, $x_{i,0}=x^*=\mathbf{1}/\sqrt{N}$,    and  $  \gamma_k={1\over k}$.   Suppose  A4 holds, and, in  addition, that

\noindent B1)    for any $i \in \mathcal{V}$,  $\{A_{i,k}\}$ is almost surely bounded;

  \noindent B2) $A$ is positive definite, and  the largest eigenvalue of $A$ has unit multiplicity with the corresponding
 unit eigenvectors denoted by  $u^{(1)}$;

\noindent B3)  for any $i\in \mathcal{V}$,    $\liminf\limits_{k\rightarrow\infty }  \mathbb{E}\left[\|y_k^TN_{i,k} z_k \| \big | \mathcal{F}_k  \right]>0$ for any $\mathcal{F}_k$-measurable  $y_k$ and  $z_k $   bounded from above and from zero, where $N_{i,k} \triangleq A_{i,k}-A_i$ and
  $\mathcal{F}_k\triangleq \sigma\{x_{i,p},i\in \mathcal{V}, 0\leq  p\leq k\}$.

 \noindent Then    \eqref{centalform}  takes place and $X_{\bot, k} \xlongrightarrow [k \rightarrow \infty]{} \mathbf{0} ~a.s.$
 Assume further that with probability one $X_k$ does not converge to $\mathbf{0}$.
Then for any $i,j\in \mathcal{V}$   and  for almost all  $\omega \in \Omega$,   the following holds:
 \begin{equation}\label{PCA_cov}
 \begin{split}
  \lim_{k \rightarrow \infty} x_{i,k}(\omega)= \lim_{k \rightarrow \infty} x_{j,k}(\omega) &=u^{(1)}  .
 \end{split}
\end{equation}
 \end{cor}
 \begin{IEEEproof}
   We prove this result by  adopting the similar procedures  as that  used in the proof of   Theorem \ref{thm1}.
Firstly, we use the auxiliary sequences defined in  \eqref{aux0}\eqref{aux1}  to prove the  finiteness of truncation numbers. Then we show that the estimates are bounded and asymptotically reach consensus. Finally,  we show that the estimates  converge  to the unit eigenvector corresponding to the largest eigenvalue.

 It is clear  that  the step size $\gamma_k ={1\over k}$ satisfies A1.   By  $f_i(x)=A_ix-(x^TA_ix)x$, we see A3 holds.
Then by \eqref{ob}\eqref{pca}  we have   \begin{equation}
  \begin{split}  \varepsilon_{i,k+1}& = O_{i,k+1}-f_i(x_{i,k}) \\&=(A_{i,k}-A_i)x_{i,k}-\left(x_{i,k}^T ( A_{i,k} -A_i)x_{i,k}\right)x_{i,k} . \nonumber
  \end{split}
\end{equation}
Then from B1 we conclude that $\{ \varepsilon_{i,k+1}I_{[ \| x_{i,k} \| \leq K]}\}$ is an  MDS with  bounded second moments  for any  $K>0$.  Then  by  the convergence theorem for  MDS \cite{Yuan1997Probability} we have that
    \begin{equation}
\label{filter5}
\begin{split}
& \sum_{k=0}^{\infty} \gamma_k \varepsilon_{i,k} I_{[ \| x_{i,k} \| \leq K]}  < \infty ~a.s.
\end{split}
\end{equation}
So,  A5 b) holds for amost all sample path $\omega$ for any agent $i$.

  Denote by $J^0$   the  set of nonzero roots of  $f(x)=Ax-(x^TAx)x$ in $\mathbb{R}^l$. Take some $u_0\in J^0$.
 Then we get $Au_0=(u_0^TAu_0)u_0,$ and hence by multiplying both sides with $u_0^T$ from the left we derive $u_0^TAu_0=(u_0^TAu_0)u_0^Tu_0.$
     Since $u_0^TAu_0>0$  by the   positive definiteness of $A$,  we have that $u_0^Tu_0=1.$  Thus,   $J^0$ is on the unit sphere.
By noticing that the  nonzero roots of   $f (\cdot)$ on the unit sphere  are the unit eigenvectors of $A$, and hence
$J^0$  coincides with   the  set of all unit eigenvectors of $A$.
Following  \cite{delyon1996general}, we define  the  Lyapunov  function as follows:
\begin{equation}
   \begin{array}{lll}
v(x)=\frac{e^{\|x\|^2}}{ x^T Ax}~~\forall x \neq \mathbf{0}. \notag
  \end{array}
   \end{equation}
Then we have
\begin{equation}\label{ }
\begin{split}
  v_x(x) & =2x v(x)-\frac{e^{\|x\|^2} 2Ax}{ (x^T Ax)^2}
 = \frac{ 2v(x)}{x^TAx}\big( (x^TAx)x-Ax \big), \nonumber
\end{split}
\end{equation}
and hence
\begin{equation}\label{lyap}
\begin{split}
f^T(x)  v_x(x) =  -\frac{ 2v(x)}{x^TAx} \| f(x)\|^2 <0~~   \forall  x \notin J^0 \cup \{\mathbf{0}\} .  \end{split}
\end{equation}

We now prove that  $ \lim\limits_{k \rightarrow \infty} \sigma_k=\sigma<\infty$.
 Assume the converse:  $ \lim\limits_{k \rightarrow \infty} \sigma_k=\infty.$

By Lemma \ref{lemma4} i), $\{\widetilde{X}_k\}$  has a convergent subsequence  $\{\widetilde{X}_{n_k}\}$
with $\widetilde{X}_{n_k}=(\mathbf{1} \otimes \textrm{I}_m) x^*$.
Since  A1, A3, A4 and A5 b) hold,   Lemmas    \ref{lemma2} and \ref{lemnoise} take place.  By noticing $\bar{x}_{n_k}=\mathbf{1}/\sqrt{N}$, from  \eqref{res22} we see  that
for sufficiently small $T$ and large $k$,
all components of  $\overline{x}_{m+1}~~\forall m: n_k \leq m \leq m(n_k,T)$ are uniformly above zero.  This together with \eqref{lyap} yields \eqref{estm41} when we follow the proof procedure for Lemma \ref{lemma3}  with $v(J)$ replaced by $v(J^0)$. It is worth noting that $v(J)$ is not defined since $v(x)$ is not defined at $x=0,$ while $v(J^{(0)})$ is well defined. Then, any nonempty interval $[\delta_1, \delta_2 ]$ with    $d([\delta_1, \delta_2 ],  v(J^0))>0 $ cannot be crossed by    $\{v(\bar{x}_{n_k}), \cdots, v(\bar{x}_{m_k})\}$  infinitely many times.

Since $v(x)$ is radially  unbounded,  there exists $c_0>1$ such that  $v(x^*)<\textrm{inf }_{\parallel x \parallel =c_0} v(x)$.   By Lemma \ref{lemma4} i),
 $\{\bar{x}_{n_k}\}$   starting from $ x^{*}$  crosses  the sphere with  $\| x \| =c_0$  infinitely many times.
Since   $v(J^0)$ is nowhere dense, there exists a nonempty interval  $[\delta_1, \delta_2 ]
    \in (v(x^*),\textrm{inf }_{\parallel x \parallel =c_0} v(x))$  with  $d\left([\delta_1, \delta_2 ],  v(J^0)\right)>0 $.
   Therefore,   $[\delta_1, \delta_2 ]$ with    $d\left([\delta_1, \delta_2 ],  v(J^0)\right)>0 $   is crossed by    $\{v(\bar{x}_{n_k}), \cdots, v(\bar{x}_{m_k})\}$  infinitely many times.  This yields a contradiction,    hence the converse assumption is not true. Thus,
   \begin{equation}\label{PCAT}
 \lim\limits_{k \rightarrow \infty} \sigma_k=\sigma<\infty.
\end{equation}

  Since \eqref{PCAT} holds, the proof for  Theorem \ref{thm1} i)  is still applicable.  As a result,   $\{X_k\}$ is  a.s. bounded and
   there  exists a positive integer $k_0$ such that the compact form \eqref{centalform}  holds.
 Then by \eqref{filter5} we obtain that $ \sum_{k=0}^{\infty} \gamma_k \varepsilon_{i,k}   < \infty~~a.s.,$
  and   hence A5 holds for almost all $\omega$.
    Since the proof for   the first result in Theorem \ref{thm1} ii) still holds,
  we  obtain    $X_{\bot, k} \xlongrightarrow [k \rightarrow \infty]{} \mathbf{0}~a.s.$
   {We further show the rate of convergence of $X_{\bot, k}$. By the   boundedness of $\{ X_k\}$ and by
   B1,   we see that $  \left \{\|  \varepsilon _{k}\| \right\}_{k\geq k_0}$ is bounded a.s.
   Then similar to the proof of  \eqref{er2}, for almost all  $\omega$,  there exist positive  constants $c_1 , c_2    $
    possibly depending  on $\omega$   such that for any $k\geq k_0:$
  \begin{equation}\label{er3}
   \begin{array}{lll}
  & \parallel X_{\bot, k+1} \parallel \leq  c_1  \rho^{k+1- k_0} + c_2  \sum_{m= k_0} ^k \gamma_m \rho^{k -m}   .
  \end{array}
   \end{equation}
   Note that
\begin{equation}
   \begin{array}{lll}
&  {1\over \gamma_k} \sum_{m= k_0} ^k \gamma_m \rho^{k -m} =\sum_{m= k_0} ^k  {k\over m} \rho^{k -m}
\\& =\sum_{m= k_0} ^k \left(1+ {k-m\over m}\right) \rho^{k -m} \\&
\leq  \sum_{m= k_0} ^k  \rho^{k -m} +\sum_{m= k_0} ^k  (k-m+1)   \rho^{k -m}
\\&\leq {1\over 1-\rho} + {d\over d \rho}  \sum_{m= k_0} ^k  \rho^{k -m+1}\leq {1\over 1-\rho} +{1\over (1-\rho)^2}  , \notag
  \end{array}
   \end{equation}
  which incorporating with \eqref{er3} yields that for almost all  $\omega$:
 \begin{equation}\label{rate-dis}
\limsup_{k\rightarrow\infty } \gamma_k^{-1}X_{\bot, k} \leq {c_2(\omega)\over 1-\rho} +{c_2(\omega)\over (1-\rho)^2}  .
\end{equation}
}

   By  recalling  that with probability one $X_{k}$ does not converge to  zero,
     for sufficiently large $k_1>k_0$  we have  that
        \begin{equation}\label{nonzero}
     X_k\neq \mathbf{0}~~\forall k \geq k_1~~a.s.
     \end{equation}
   This is because  the converse assumption
   $X_{k_2}= \mathbf{0} $ for some      $k_2>k_1 $ implies that   $X_k= \mathbf{0}~\forall k \geq k_2$  by \eqref{pca}  and \eqref{notr}.  This contradicts with the assumption that $X_{k}$ does not converge to zero. So, \eqref{nonzero} holds. Hence
   by  $X_{\bot, k} \xlongrightarrow [k \rightarrow \infty]{} \mathbf{0}$ it follows that  $x_k$ does not converge to zero.
{ Since A1, A3, A4, A5 and  \eqref{lyap} hold,     similar  to the proof  of   Theorem \ref{thm1} iii) we have  that
    for almost all $\omega \in \Omega:$
 \begin{equation}\label{converge-connected}
  \lim\limits_{k \rightarrow \infty} d(x_k(\omega),J^*)=0,
\end{equation}
where  $J^{*} $ is a connected  subset $ J^0$.

Similar to the procedure for deriving \eqref{centrelized2},  we can rewrite \eqref{centalform}   in the  following  centralized form:
\begin{equation} \label{centr-from}
  \begin{split}
x_{k+1}=x_k +  { \gamma_k \over N } \left( f(x_k)+\zeta_{k+1}+   \beta_{k+1}  \right),
    \end{split}
\end{equation}
 where  $\beta_{k+1}\triangleq \sum_{i=1}^N  \left(  f_i(x_{i,k})- f_i(x_k) \right) $ and
$  \zeta_{k+1}\triangleq  (\mathbf{1}^T \otimes \mathbf{I}_l)   \varepsilon_{k+1}  .$}
  Note that $\{\zeta_k\}$ is an MDS with bounded second moments by B1, $\beta_k=O(\gamma_k)$ a.s. by (64), and     $\{x_k\}$ is bounded.
Denote by $ \lambda^{(1)},\cdots,  \lambda^{(l)} $ the eigenvalues of $A$ in a nonincreasing order. Then
 $ \lambda^{(1)}>  \lambda^{(p)}~\forall p>1$  by B2.
Multiplying both sides of   \eqref{centr-from} from the left with $(u^{(1)})^T$,
by  noticing  that $ (u^{(1)})^T A= \lambda^{(1)} (u^{(1)})^T  $  and  $f(x)=Ax-(x^TAx)x$, we  obtain the following:
\begin{align}
(u^{(1)})^T x_{k+1}&=(u^{(1)})^T x_k +  { \gamma_k \over N }   \left( \lambda^{(1)}   -x_k^TAx_k  \right)(u^{(1)})^T  x_k \notag
\\&+ { \gamma_k \over N } (u^{(1)})^T  \left(\zeta_{k+1}+   \beta_{k+1}  \right).  \label{bd-centr-from}
    \end{align}
Notice that $u^{(1)}$ is orthogonal to other eigenvectors of $A$ and that $x_k$ converges to some  unit eigenvector of A by \eqref{converge-connected}. So, to prove \eqref{PCA_cov} it suffices
  to show that $\mathbb{P}(\Omega_1)=0$, where  $\Omega_1=\{\omega: (u^{(1)})^T x_k(\omega) \to 0\}.$
 By definition of  $ \varepsilon_{i,k+1}$,  we get
    \begin{equation}
  \begin{split}
  &\mathbb{E}\left[\|(u^{(1)})^T  \varepsilon_{i,k+1} \| \big | \mathcal{F}_k  \right]
  \\&=\mathbb{E}\left[\|(u^{(1)})^TN_{i,k}x_{i,k}-\left(x_{i,k}^T N_{i,k}x_{i,k}\right) (u^{(1)})^Tx_{i,k}\| \big | \mathcal{F}_k  \right]
  \\&\geq \mathbb{E}\left[\|(u^{(1)})^TN_{i,k}x_{i,k}\| \big | \mathcal{F}_k  \right]
 \\& - \left \|(u^{(1)})^Tx_{i,k} \right\| \cdot \mathbb{E}\left[\| \left(x_{i,k}^T N_{i,k}x_{i,k}\right)  \| \big | \mathcal{F}_k  \right]   . \nonumber
  \end{split}
\end{equation}
Then  by B3,  $\liminf_k \mathbb{E}\left[\|(u^{(1)})^T  \varepsilon_{i,k+1} \| \big | \mathcal{F}_k  \right]  >0 $ on $\Omega_1$,
and hence $\liminf_k \mathbb{E}\left[\|(u^{(1)})^T  \zeta_{k+1}   \| \big | \mathcal{F}_k  \right]  >0 $ on $\Omega_1$.
Note that on $\Omega_1$,  by \eqref{converge-connected}   we have that  $\lambda^{(1)} -x_k^TAx_k \to \lambda^{(1)}-\lambda^{(p)}>0 $ for some $p>1.$
Then by \eqref{bd-centr-from}  and  \cite[Lemma 1]{chen2011recursive}, we obtain  that  $\mathbb{P}(\Omega_1)=0$.  The  proof is completed.
 \end{IEEEproof}

\begin{rem} \cite[Proposition]{chen2011recursive}
 For any $i\in \mathcal{V},$  assume that $(N_{i,k},\mathcal{F}_k)$ with $ N_{i,k}= \left \{N_{i,k}(p,q) \right \}_{p,q=1}^l$ is an MDS,  $\mathbb{E}\left[N_{i,k+1}(p,q)N_{i,k+1}(s,t)\big| \mathcal{F}_k\right]=0$
whenever $(p,q) \neq (s,t)$, and $\liminf_k \mathbb{E}\left[\| N_{i,k+1}(p,q)\|^2\big| \mathcal{F}_k\right]   >0~\forall p,q=1,\cdots,l.$ Then B3 holds.
\end{rem}

In  Corollary \ref{thm-PCA}   ``non-convergence of $X_k$ to zero" is assumed rather than proved, but the numerical simulations substantiate this assertion. This is not surprising because if we compare \eqref{centr-from} with $u_k^{(i)}$ given in \cite{Chen_Zhao}, the estimates for unit eigenvectors corresponding to eigenvalues arranged in the decreasing order of a noisy observed matrix $A$  (pp. 289--316 in \cite{Chen_Zhao}), we find that \eqref{centr-from} and the recursive expression ((5.1.9) of \cite{Chen_Zhao}) of $u_k^{(1)}$ are in a complete similarity while the latter converges to the unit eigenvector corresponding to the largest eigenvalue of $A$ as proved in \cite{Chen_Zhao}. It is worth noting that normalization technique used in \cite{Chen_Zhao} automatically excludes estimates from tending to zero, but this can also be achieved by truncations as done in  \cite[Section 6.2]{Chen_2002}  when solving the adaptive stabilization problem, where not only the upper side truncations are used to prevent estimates from diverging to infinity but also the lower side truncations for estimates not tending to zero. We might apply this kind of idea to the  formulated distributed  root-seeking problem, then the modified  distributed  stochastic approximation algorithm would converge to nonzero roots of sum  functions.

  \begin{exa} Distributed Principal Component Analysis

 Let  $N=1000$.   The   matrix  $W(k)$ is given by the following:
 \begin{equation} \begin{array}{lll}
&W(3k-2)= \begin{pmatrix}
    W_1  &   \mathbf{0} \\
     \mathbf{0} &  \textrm{I}_{N/2}
\end{pmatrix},~ W(3k-1)= \begin{pmatrix}
      \textrm{I}_{N/2} &   \mathbf{0} \\
     \mathbf{0} &W_2
\end{pmatrix},  \\
&W(3k)= \begin{pmatrix}
    \frac{1}{2}  \textrm{I}_{N/2} &  \frac{1}{2} \textrm{I}_{N/2}\\
   \frac{1}{2}\textrm{I}_{N/2}&  \frac{1}{2}\textrm{I}_{N/2}
\end{pmatrix},  \nonumber
\end{array}\end{equation}
where matrices $W_1\in \mathbb{R}^{N/2\times N/2}$ and  $W_2\in \mathbb{R}^{N/2\times N/2}$ are  doubly stochastic.     Further, let   $W_1$ and $W_2$ be    the  adjacency  matrices of some strongly connected digraphs.
Thus, A4 holds.

 Each sensor  $i=1,2, \cdots, N$ has access to a  $9$-dimensional i.i.d.
 Gaussian sequence $\{u_{i,k}\}$ with zero mean.
   Set $\gamma_k=\frac{1}{k}$, $M_k=2^k$.  Let the   sequence $\{x_{i,k}\}$
   be produced by \eqref{comp1}--\eqref{step2} with $O_{i,k+1}$ defined by \eqref{pca} and with
   $x_{i,0}= x^{*}=\mathbf{1}/\sqrt{N}$.
    Denote by $x_k^i$ the   $i$th component  of   $x_k=\frac{1}{N} \sum_{i=1}^N x_{i,k}$,    and by   $e(k)=\sum_{i=1}^N \| x_{i,k}-u^{(1)} \|_2/N$ the average of   2-norm errors  for all agents  at time $k$.
 The    estimates  $x_k$ are demonstrated  in  Fig. \ref{Figure A1}   with the solid  lines denoting  the true values and   the dotted, dashed, and dash-dot lines  representing  the estimates.  Meanwhile, the  estimation errors  $e(k)$ are demonstrated in Fig. \ref{Figure A1}.  From the figure it is seen that the estimates converge to the unit eigenvector corresponding to the largest eigenvalue of $A.$

\begin{figure}[!htb]
    \centering
  \includegraphics[width=3.5in]{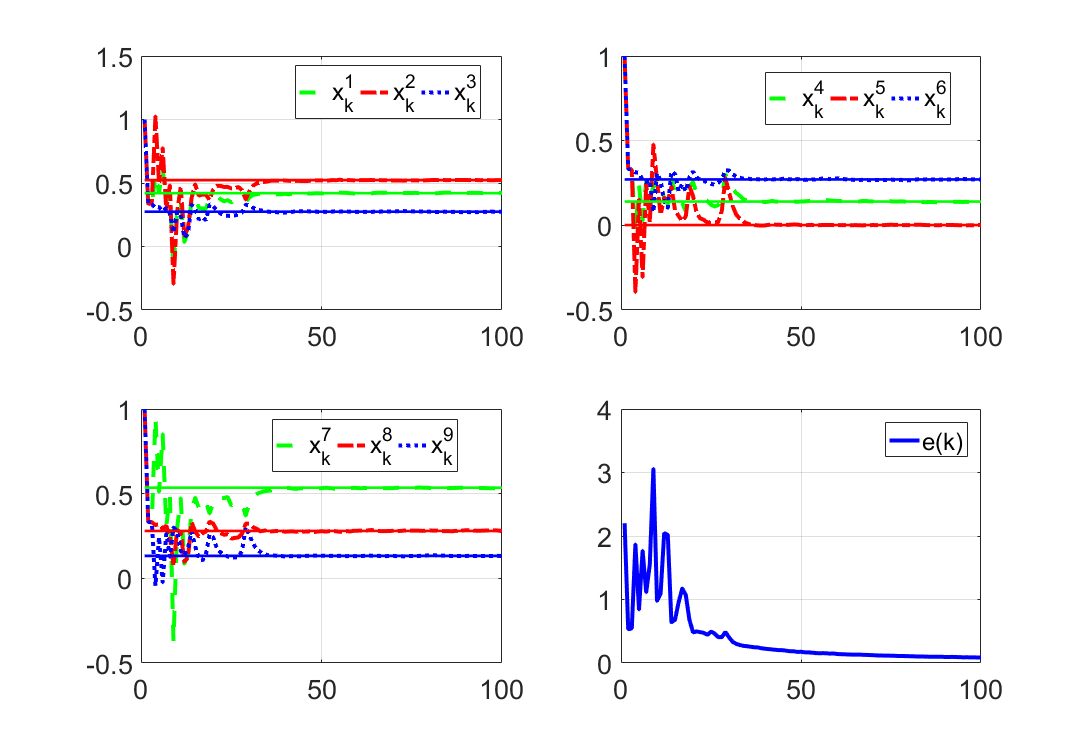}
  \caption{The   estimation sequences   and estimation errors} \label{Figure A1}
\end{figure}
\end{exa}

  For a comparison, we also present the simulation results computed by DSAA  proposed in \cite{DSA} for the same example.
The  initial value for each agent,     $  \alpha_k$ and $ \gamma_k  $
 are set to be the same as  those used in DSAAWET.     The  simulation results are given in Table II
 with   $ e+11$ denoting the exponential  $10^{11}$,  from which it is seen that the estimates are unbounded.  When    a smaller    step size, say $\gamma_k=\frac{ 1}{ k+20 }$,  is taken,  then the  estimates are bounded  for some sample paths and
  are unbounded  for  others.
 It is hard to indicate how small the  step-sizes should be for  a given sample  path in advance.
As a matter of fact, DSAAWET   adaptively   chooses   step-sizes for any sample path.

\renewcommand{\arraystretch}{1.1}
 \begin{table}[!htb]
\centering
     \begin{tabular}{|c|c|c|c|c|c|c|}
        \hline
     $k$   & 0 & 1 & 2& 3& 4 \\ \hline
     $x_k^1$ & 1/3& 1.9 & -587 &  $5.46e+{11}$ & $-1.49e+{41}$  \\ \hline  
  $x_k^2$ & 1/3& 3.2 & -593 &  $3.64e+{11}$ & $-9.89e+{40}$  \\ \hline  
  $x_k^3$ & 1/3& 0.277 & -45.4 &  $-9.9e+{10}$ & $2.93e+{40}$  \\ \hline
  $x_k^4$ & 1/3& -1.26 & 402 &  $-6e+{11}$ & $1.72e+{41}$  \\ \hline
  $x_k^5$ & 1/3& -3.26 & 736 &  $-2.23e+{11}$ & $4.76e+{40}$  \\ \hline
  $x_k^6$ & 1/3& 0.25 & -17.8 &  $7.86e+{10}$ & $-2.23e+{40}$  \\ \hline
  $x_k^7$ & 1/3& 3.3 & -711 &  $5.02e+{11}$ & $-1.34e+{41}$  \\ \hline
  $x_k^8$ & 1/3& 0.19 & 2.48 &  $-1.76e+{11}$ & $5.54e+{40}$  \\ \hline
  $x_k^9$ & 1/3& -1.63 & 397 &  $-2.37e+{11}$ & $6.38e+{40}$  \\ \hline
    \end{tabular}
    \caption{Estimates  produced   by  DSAA  in \cite{DSA}  }
\end{table}
\subsection{Distributed Gradient-free Optimization}

We now  consider  the convergence of DSAAWET applied to the  distributed gradient-free  optimization
problem.

\begin{cor}\label{thmopt}
Let $ \{ x_{i,k}\} $ be produced by \eqref{comp1}-\eqref{step2} with  $O_{i,k+1}$ given by  \eqref{opt2}
for any initial value  $   x_{i,0} $.
 Assume   that     A4 holds, and, in addition,   \\
 C1 i)  $\gamma_k>0, \sum_{k=1}^{\infty } \gamma_k =\infty$, and   $ \sum_{k=1}^{\infty } \gamma_k^p <\infty $
 for some   $p \in (1,2]$;

 ~ii) $\alpha_k>0$, and $\alpha_k\rightarrow 0$ as $k \rightarrow  \infty$;

 ~iii) $ \sum_{k=1}^{\infty }  \gamma_k^2 / \alpha_k^2  <\infty$;\\
 C2 i) $  c_i(\cdot), i=1,2,\cdots,N$ are continuously  differentiable, and  $\nabla c(\cdot)$ is locally Lipschitz continuous;

 ~~ii) $c(\cdot)$ is convex with a  unique global  minimum  $x_{min}$;

      ~iii)   $c(x^{*}) <  \sup_{ \| x \| =c_0} c(x) $ and   $\| x^{* } \| <c_0$   for some positive constant $c_0$,
      where   $x^{*}$  is used in \eqref{comp2} \eqref{step1};\\
 C3  $\xi_{i,k+1}= \xi_{i,k+1}^{+}- \xi_{i,k+1}^{-}  $  is independent of $ \{ \Delta_{i,s}, s=1,2,\cdots,k\}$ for any
 $k\geq 0$, and $\xi_{i,k+1}$ satisfies one  of the following  two conditions:

 ~i) $\sup_k | \xi_{i,k } | \leq \xi_i ~ a.s.,$ where $\xi_i$ may be random;

 ~ii) $\sup_k  E [ \xi_{i,k+1}^2] < \infty$  \\
  for any $i=1,2,\cdots N$. Then
  \begin{equation} \label{opt6}    x_{i,k} \xlongrightarrow [k \rightarrow \infty] {} x_{min} \quad a.s.
  \end{equation}
  for any $i=1,2,\cdots, N$.
\end{cor}
\begin{IEEEproof}
  To apply    Theorem \ref{thm1}   to this problem, we have to verify  conditions A1-A5.

By C1  the step-size $\{ \gamma_k\}$ satisfies A1.
   Since $c(\cdot)$ is convex and  differentiable,
   $x_{min}$ is the global   minimum  of $c(\cdot) $ if and only if
   $ \nabla c(x_{min})=0.$
   So, the original   problem  \eqref{opt1}  is equivalent to  finding the root $J=\{ x_{min}\}$ of
   $ f(\cdot)=\sum_{i=1}^N f_i(\cdot)/N$ with $f_i(\cdot) =- \nabla c_i(\cdot)$.
  By setting   $v(\cdot)=c(\cdot)$, we  obtain that  $f^T(x)v_x(x)=- \| \nabla c(x) \|^2 /N $ and  $v(J)=\{ c(x_{min})\}$.
  So,   A2  a) and A2 b) hold.  By C2  iii) we have A2  c).  It is clear that  C2 i)  and iii) imply A3.

  By  \eqref{ob}    we  have that
  \begin{equation}\label{DO_noise}
\varepsilon_{i,k+1}=O_{i,k+1}+\nabla c_i(x_{i,k}) ,
\end{equation}
where  $O_{i,k+1}$ is   given by    \eqref{opt2}.
 The analysis for the observation noise $\{\varepsilon_{i,k+1}\}$ is the same as that  given    in  \cite{Chen_1999},
  and hence  is omitted  here.  It is shown in  \cite{Chen_1999}    that
  $\{ \varepsilon_{i,k+1} \}~\forall  i \in \mathcal{V}$  is not an  MDS,   while using the  local Lipschitz continuity  of $\nabla c(\cdot)$, C1  and  C3, by  \cite[Lemma 2]{Chen_1999},  the following limit takes place:
\begin{equation}  \label{opt5}
\begin{split}
  & \lim_{T \rightarrow 0} \limsup_{k \rightarrow \infty }\frac{1}{T}
   \left \|   \sum_{n= k}^{m( k,t_k)}\gamma_{m } \varepsilon_{i,n+1} I_{ [\parallel x_{i,n} \parallel \leq K ]} \right \| =0
 \\& ~~~  \forall t_k \in [0,T]   \textrm{   for any positive integer }K , ~~a.s.
  \end{split}
\end{equation}
Therefore,   A5  b) holds  with probability one   for any $i \in \mathcal{V}$.  By Theorem \ref{thm1} i), we see that
  $\{ x_{i,k}\}$  is bounded almost surely.  So,
A5 a) holds   by taking $t_k= \gamma_k$ in  \eqref{opt5}.

In summary,  we have shown that  A1-A5 hold. Since $J=\{ x_{min}\}$,  by  Theorem \ref{thm1} ii)   we derive \eqref{opt6}.
    \end{IEEEproof}

\begin{rem}
If the convex function $c(\cdot)$ is allowed to have non-unique  minima while
 other conditions in  Corollary    \ref{thmopt} remain  unchanged,  then by
Theorem \ref{thm1} we know that   $x_{i,k} ~\forall i \in \mathcal{V}$ converge to a connected subset $J^{*} \subset J
\triangleq \{ x \in \mathbb{R}^l: \nabla c(x)=0\}$.
\end{rem}

\begin{rem}
Compared with  the biased mean error-bounds for the  local weighted average vector given  in \cite{Yuan_NN_15}\cite{li2015gradient},
 the almost sure convergence of the iterates to the global optimal solution  set is  established here.
Regarding the problem  settings: i) each  nonsmooth  objective  is
assumed to be convex and Lipschitz continuous   in  \cite{Yuan_NN_15}\cite{li2015gradient},
  while in this work,  each objective function is not necessarily   convex but
required to be  continuously  differentiable with the gradient function being locally Lipschitz continuous;
ii)  a convex set constraint optimization is considered  in  \cite{Yuan_NN_15}\cite{li2015gradient}
while   the unconstrained optimization  is    considered  here;
iii)  the noisy observations of the objective functions  are considered  here in contrast to    the exact observations
discussed   in \cite{Yuan_NN_15}\cite{li2015gradient}.
\end{rem}

\begin{exa} Distributed Gradient-free Optimization

 Consider the  network of three agents with  local cost  functions given by
\begin{equation}\label{sf}
\begin{split}
& L_1(x,y)=x^2+y^2+10\sin(x) ; \\
 & L_2(x,y)=(x-4)^2+(y-1)^2-10\sin (x);\\
 & L_3(x,y)=0.01(x-2)^4+(y-2)^2. \nonumber
\end{split}
\end{equation}

Let  the communication relationship   among the agents  be described by a
  strongly connected digraph with the adjacency matrix being doubly stochastic.
 The task of the network is  to find the minimum   $(x^0,y^0)=(2,1)$ of the cost function $  L(x,y) =\sum_{i=1}^3 L_i(x,y)$.
Though each local cost function is non-convex, the    global  cost function $L(x,y)$  is convex.

Let the observation noise for the cost function of each agent be a sequence of i.i.d. random vectors  $\in {\cal N}(0,I).$
The first and second  component of the initial values for all  agents are set to be mutually independent and uniformly distributed over the intervals  $[-2,6]$  and  $[-2,4]$, respectively.
Set $x^*= (-1,4)^T$, $\gamma_k=\frac{2}{k}$, $\alpha_k= \frac{1}{k^{0.2}}$,  and $M_k=2^k$.
Let the  estimates for the minimum
   be produced by \eqref{comp1}--\eqref{step2} with $O_{i,k+1}$ defined by \eqref{opt2}.
 The estimates of $x^0$  and  $y^0$ produced by  the three agents are demonstrated in  Fig. \ref{Figure A2},
 where   $x_i(k) $ and $ y_i(k) $ denote    the estimates for $x^0$ and $y^0$, respectively,
 given at agent $i$ and  at time $k$. From the figure  it is seen that the estimates given by all agents converge to the minimum, which is consistent with  the theoretic result.

\begin{figure}[!htb]
     \centering
  \includegraphics[width=2.8in]{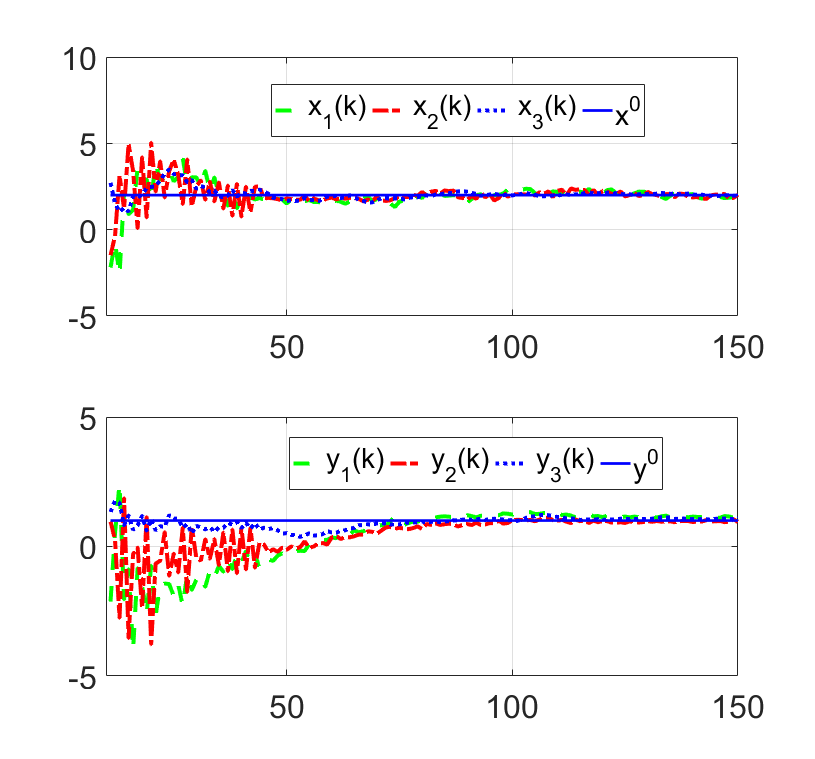}
   \caption{ Trajectories of   the estimates } \label{Figure A2}
\end{figure}
\end{exa}

\section{ Concluding Remarks}
In this paper,   DSAAWET is defined to solve the  formulated distributed root-seeking problem.
   The estimates are shown  to converge to a consensus set   belonging  to a connected subset of   the root set.  Two  problems as examples of those which   can be solved by DSAAWET are demonstrated with  numerical  simulations  provided.

For further research it is of interest to analyze the convergence rate of the proposed algorithm, and to  consider
the convergence properties of DSAAWET  over random networks taking into account the possible packet loss in communication.

\appendices
\numberwithin{equation}{section}

\section{Proof Lemma \ref{lemau2}}\label{PLA}

Before proving the lemma,  we show that 
for  any $ k \in [ \tau_m , \tau_{m+1})$, the following assertions hold:
 \begin{flalign}
&~~~\textrm{ i)}  ~~~ \tilde{x}_{i,k}=x^* , \quad   \tilde{\varepsilon}_{i,k+1}=-f_i(x^*)\textrm{ if }\sigma_{i,k}<m;  &  \label{case1}\\
&~~~\textrm{ ii)}  ~~~ \tilde{x}_{i,k}=x_{i,k} , \quad   \tilde{\varepsilon}_{i,k+1}=\varepsilon_{i,k+1}~ \textrm{ if }\sigma_{i,k}=m;  &  \label{case2} \\
&~~~\textrm{ iii)}  ~~~ \tilde{x}_{j,k}=x^*   \textrm{ if }\sigma_{j,k-1}<m; &  \label{case3}\\
&~~~\textrm{ iv)}  ~~~ \tilde{x}_{j, k+1}=  x^{*}~~ \forall j \in \mathcal{V} ~~\textrm{ if } \sigma_{k+1}=m+1.   &  \label{sub2}
\end{flalign} 
{\bf Proof:}
i) From  $\sigma_{i,k} <m$ we see that the truncation number of agent  $i$
is smaller than $m$   at time $k$, then  by the definition of $\tau_{i,m}$ we derive  $\tau_{i,m} >k$.
   Thus,  $ \tilde{\tau}_{i,m}= \tau_{i,m} \wedge \tau_{m+1} >k,$ and  hence  from  \eqref{aux0} we
   conclude \eqref{case1}.

ii) From   $\sigma_{i,k}=m$  by definition we have $\tau_{i,m} \leq k$,
 and hence $\tilde{\tau}_{i,m} = \tau_{m+1} \wedge \tau_{i,m } =\tau_{i,m}  \leq  k. $
Then by   \eqref{aux1}  it is clear  that \eqref{case2} holds.

iii) By  $\tau_m \leq  k <  \tau_{m+1}$ we see   $\sigma_{j,k} \leq m ~ \forall  j \in \mathcal{V}$.
 We show \eqref{case3} separately for the cases   $\sigma_{j,k} = m$ and  $\sigma_{j,k} < m$.
 1)  Let us first consider the case $\sigma_{j,k} = m$. Since  $\sigma_{j,k-1} < m   $ and $\sigma_{j,k} = m$,
        by \eqref{sub1} we obtain $x_{j,k}=x^{*}$.
Hence  from \eqref{case2} we see  $\tilde{x}_{j,p} =x_{j,p}=x^{*}$.
2)   We now consider the case $\sigma_{j,k} < m$.   By  \eqref{case1} we see $\tilde{x}_{j,k}=x^{*}$, which is
the assertion of  \eqref{case3}.

iv) From  $ k \in [ \tau_m , \tau_{m+1})$    we see $\sigma_k=m$.
 Thus from $ \sigma_{k+1}= m+1$ by  definition
   we derive  $\tau_{m+1} =k+1$, and hence $  k+1 \in [ \tau_{m+1} , \tau_{m+2})$.
 By $\sigma_k=m$  we see  $\sigma_{j,k} < m+1 ~ \forall  j \in \mathcal{V}$,
  then we derive  \eqref{sub2}   by \eqref{case3}.
\hfill $\blacksquare$

 {\bf Proof of Lemma \ref{lemau2}:}
We prove the  lemma by induction.

We first prove    \eqref{algorithm01}-\eqref{algorithm03}  for $k=0$.
Since $0 \in [\tau_0, \tau_1)\textrm{ and } \sigma_{i,0}=0 ~ \forall i \in \mathcal{V}$,
by \eqref{case2} we  derive
 $ \tilde{x}_{i,0} =  x_{i,0}, ~  \tilde{\varepsilon}_{i, 1} = \varepsilon_{i,1}   ~ \forall i \in \mathcal{V}$.
Then by noticing  $\hat{\sigma}_{i,0}=\sigma_{i,0}=0 ~ \forall i \in \mathcal{V}$,
from \eqref{comp2}   \eqref{algorithm01} we see
\begin{equation} \label{equalityhat}
\hat{x}_{i, 1}= x'_{i,1}~~ \forall i \in \mathcal{V}.
\end{equation}
 We  now show  that $  \tilde{x}_{i,1}$ and $\sigma_1$ generated by  \eqref{algorithm01}-\eqref{algorithm03}
 are consistent with their definitions \eqref{aux0} \eqref{aux1} \eqref{largest} by considering the following two cases:

  i) There is no truncation  at $k=1$, i.e.,  $\sigma_{i,1}=0  ~\forall i \in \mathcal{V}$.
   In this case, from $\hat{\sigma}_{i,0} =0 ~ \forall i \in \mathcal{V}$  by \eqref{step2} we have  $  \| x'_{i,1} \| \leq M_0 ~\forall i \in \mathcal{V}$, and hence      $ \| \hat{x}_{i, 1}\| \leq M_0 ~\forall i \in \mathcal{V}$   by \eqref{equalityhat}.
Then  $ x_{i, 1}= x'_{i,1} $  by \eqref{step1},  $  \tilde{x}_{i, 1} = \hat{x}_{i, 1}$  and $\sigma_1=0$
      by \eqref{algorithm02}   and  \eqref{algorithm03}, respectively.
       These together with \eqref{equalityhat} imply  that $  \tilde{x}_{i,1}= x_{i, 1} ~ \forall i \in \mathcal{V},$ which is consistent
with the definition of    $\tilde{x}_{i,1} $  given by \eqref{aux1} since   $\tilde{\tau}_{i,0} \leq  1 < \tau_1$.
By   \eqref{largest} we see $\sigma_1=\max_{i \in \mathcal{V}} \sigma_{i,1}=0,$ which is consistent
with that derived from    \eqref{algorithm01}-\eqref{algorithm03}.

  ii) There is a  truncation   at $k=1$ for some agent $i_0$, i.e.,  $\sigma_{i_0,1}=1.$
 In this case,  by \eqref{step1}  \eqref{step2} we derive
  $ x_{i_0,1}=x^{*},  ~ \| x'_{i_0,1} \|>M_0$,
  and hence $\|  \hat{x}_{i_0, 1} \|  > M_0$ by  \eqref{equalityhat}.
 Therefore,   $   ~ \tilde{x}_{i, 1} = x^{*} ~\forall i \in \mathcal{V}
\textrm{ and }\sigma_1=1$ by    \eqref{algorithm02}    \eqref{algorithm03}.
By    \eqref{largest}   from   $\sigma_{i_0,1}=1 $ we derive   $\sigma_1=\max_{i \in \mathcal{V}} \sigma_{i,1}=1$.
    Since $ 0 \in [ \tau_0 , \tau_{ 1})$ and  $\sigma_1=1$,
   by \eqref{sub2} we see $\tilde{x}_{i, 1} = x^{*}  ~ \forall i \in \mathcal{V}$.
   Thus,  $  \tilde{x}_{i,1}$ and $\sigma_1$  defined by \eqref{aux0} \eqref{aux1} and  \eqref{largest}  are  consistent with    those
     produced  by   \eqref{algorithm01}-\eqref{algorithm03}.

 In summary,    we have proved the lemma  for $k=0$.

 Inductively,  we assume   \eqref{algorithm01}-\eqref{algorithm03} hold  for $k=0,1,\cdots,p$.
 At the  fixed  sample path $\omega$ for a given  integer  $  p  $  there exists a unique integer   $m $
  such that $\tau_m \leq  p <  \tau_{m+1}$.
We  intend to show   \eqref{algorithm01}-\eqref{algorithm03} also  hold  for
   $k=p+1$.

Before doing this, we  first express  $  \hat{x}_{i,p+1}~\forall i \in \mathcal{V}$  produced by \eqref{algorithm01}
for  the following two cases:

Case 1:    $\sigma_{i,p} <m$.   Since $p \in [\tau_m, \tau_{m+1}), $ by \eqref{case1} we see
\begin{equation}\label{lem331}
  \tilde{x}_{i,p}=x^* , \quad   \tilde{\varepsilon}_{i,p+1}=-f_i(x^*).
\end{equation}
From  $\sigma_{i,p} <m$ it follows that  $\sigma_{j,p-1} < m ~\forall j \in N_i(p)$,
because otherwise, there would exist a $j_1 \in N_i(p)$
 such that  $\sigma_{j_1,p-1} \geq m$. Hence by  \eqref{comp1} \eqref{step2}
we would  derive $\sigma_{i,p} \geq \hat{\sigma}_{i,p-1}  \geq  \sigma_{j_1,p-1} \geq  m  $,
yielding  a contradiction.

From  $\sigma_{j,p-1} < m ~\forall j \in N_i(p)$  by \eqref{case3} we have
 $\tilde{x}_{j,p}=   x^*~ ~\forall j \in N_i(p) ,
$ which incorporating with   \eqref{algorithm01} \eqref{lem331} yields
 \begin{equation}\label{A11}
    \hat{x} _{i,p+1}=  x^*  \quad  \forall  i: \sigma_{i,p}<m .
 \end{equation}

 Case 2:   $ \sigma_{i,p}= m$. By  $\tau_m \leq  p <  \tau_{m+1}$ we see  $\sigma_{j,p} \leq m ~ \forall  j \in \mathcal{V}$,
   and hence by  \eqref{comp1}  we  obtain
 \begin{equation}\label{hatd}
    \hat{\sigma}_{i,p}=m~ ~ \forall  i: \sigma_{i,p}=m .
 \end{equation}
Then by  \eqref{comp2}
 \begin{equation}\label{aux2}
 \begin{array}{lll}
    x'_{i,p+1} &=& \sum_{ j \in  N_{i }(p)}  \omega_{ij}(p) ( x_{j,p} I_{[ \sigma_{j,p} = m]} +
      x^* I_{[ \sigma_{j,p} < m]}) \\&  +& \gamma_{ p } (f_i( x_{i,p})+\varepsilon_{i,p+1}).
      \end{array}
 \end{equation}

 From  $\sigma_{i,p} =m$ and  $p \in [\tau_m, \tau_{m+1}), $    by   \eqref{case2}  it is clear  that
 \begin{equation}\label{express1}
    \tilde{x}_{i,p}=x_{i,p}, ~ ~ \tilde{\varepsilon}_{i,p+1}=\varepsilon_{i,p+1} .
 \end{equation}
By the first term in \eqref{case1} and \eqref{case2},    for any $j \in N_i(p)$ we have
  \begin{equation}\label{express2}
    \tilde{x}_{j,p}=x_{j,p}~\textrm{ if }\sigma_{j,p} = m, ~ ~   \tilde{x}_{j,p}=x^{*} ~\textrm{ if }\sigma_{j,p} < m.
 \end{equation}
 Substituting \eqref{express1} \eqref{express2} into  \eqref{algorithm01},   from \eqref{aux2} we  derive
     \begin{equation}\label{eq111}
 \begin{split}
   \hat{x}_{i,p+1}=  x'_{i,p+1}  \quad   \forall  i: \sigma_{i,p}=m .
      \end{split}
 \end{equation}

Since $\tau_m \leq  p <  \tau_{m+1}$,  from $p <  \tau_{m+1}$  it follows that
$\sigma_p < m+1$, and hence $\sigma_{p+1} \leq  m+1$,
while from $\tau_m \leq  p $ it follows that
$\sigma_p=m$, and hence $\sigma_{p+1}\geq m$. Thus, we have $m\leq \sigma_{p+1}\leq  m+1$.

We  now show that   $  \tilde{x}_{i, p+1}$ and $\sigma_{p+1}$ generated by  \eqref{algorithm01}-\eqref{algorithm03}
 are consistent with their definitions  \eqref{aux0} \eqref{aux1} \eqref{largest}.
  We prove this  separately  for the   cases  $\sigma_{p+1}=m+1$ and  $\sigma_{p+1}=m$.

Case 1:    $\sigma_{p+1}=m+1$.
  We  first  show
 \begin{equation}\label{if}
   \sigma_{i,p+1} \leq m  ~\textrm{ if }\sigma_{i,p} <m
 \end{equation}
separately for  the  following  two cases 1) and 2):   1)   $\sigma_{i,p} <m$   and  $ \sigma_{j,p} <m ~\forall j \in N_i(p)$. For  this case
 by \eqref{comp1} we derive   $\hat{\sigma}_{i,p} <m$, and hence    $\sigma_{i,p+1} \leq  \hat{\sigma}_{i,p}+1 \leq  m$  by  \eqref{step2}.
2)    $\sigma_{i ,p} <m$   and  $ \sigma_{j,p} =m ~\textrm{ for some }j \in N_i(p)$.
 For  this case we obtain   $\hat{\sigma}_{i,p} =m, ~x'_{i,p+1}=x^{*}$  by \eqref{comp1} \eqref{comp2}.  Since
 $   \| x^{*}\|  \leq  M_0 \leq M_m $,   then  by \eqref{step2}   $\sigma_{i,p+1} =\hat{\sigma}_{i,p} =m$.
 Thus, $\sigma_{i,p+1}\leq m $ when $\sigma_{i,p} <m$.  Thereby  \eqref{if} holds. This means that
  \begin{equation}\label{onlyif}
     \sigma_{i,p+1} =m+1\textrm{ only if } \sigma_{i,p} =m.
\end{equation}

By definition from $\sigma_{p+1}=m+1$  we know that
  there exists  some $  i_0 \in\mathcal{V}  $ such that   $\sigma_{i_0, p+1}=m+1$.
Then  $\sigma_{i_0, p }=m$ by \eqref{onlyif}, and hence   $\hat{\sigma}_{i_0, p }=m$ from \eqref{hatd}.
 Then  from $\sigma_{i_0, p+1}=m+1$   by \eqref{step2} we derive  $\parallel  x'_{i_0, p+1} \parallel > M_m$,
 and hence
  $  \parallel  \hat{x}_{i_0, p+1} \parallel =\parallel  x'_{i_0, p+ 1} \parallel  >M_m $  by  \eqref{eq111}.
Then from \eqref{algorithm02} \eqref{algorithm03}   we derive   $   \tilde{x} _{i, p+1}=x^*~~ \forall i \in \mathcal{V} $
and  $\sigma_{p+1}=m+1 $, which is consistent
with    $\sigma_{p+1}$ defined by \eqref{largest}. Since  $\sigma_{ p+1}=m+1$ and   $p \in [\tau_m, \tau_{m+1}), $  by \eqref{sub2} or  from \eqref{aux0} \eqref{aux1}, we see  $   \tilde{x} _{i, p+1}=x^*~  \forall i \in \mathcal{V} $.
This is consistent with    that   produced  by   \eqref{algorithm01}-\eqref{algorithm03}.

Case 2:   We now consider the case $\sigma_{p+1}=m$.   In this case,   $\sigma_{i  , p+1} \leq m ~ \forall i \in \mathcal{V}.$
 By  \eqref{hatd},  from  \eqref{step1}  \eqref{step2} we see
   \begin{equation}\label{B20}
\parallel  x'_{i , p+1} \parallel \leq  M_m ,   ~x_{i,p+1}= x'_{i, p+1} ~~ \forall i :\sigma_{i, p }=m .
   \end{equation}
 So, by  \eqref{eq111}  we have
 \begin{equation}\label{B23}
   \| \hat{x}_{i,p+1} \| =  \| x'_{i,p+1}\| \leq M_m  \quad  \forall  i:\sigma_{i, p }=m .
 \end{equation}
 From  $\parallel x^* \parallel  \leq M_0 \leq  M_m$  and   \eqref{A11} we derive
$  \parallel \hat{x}_{i,p+1} \parallel \leq  M_m ~  \forall  i:\sigma_{i, p }<m ,$
which incorporating  with   \eqref{B23}  yields
 $\parallel  \hat{x}_{i ,p+ 1} \parallel  \leq M_m
  ~  \forall i \in \mathcal{V}.$ Then from \eqref{algorithm02}   we have
     \begin{equation}\label{B21}
        \tilde{x} _{i,p+1}= \hat{x}_{i,p+1}  \quad \forall i \in \mathcal{V}, ~~ \sigma_{p+1}=m.
     \end{equation}
Thus  $\sigma_{p+1} $  is  consistent with that  defined by \eqref{largest}.

It remains to show that  $  \tilde{x}_{i,p+1}$  generated by  \eqref{algorithm01}-\eqref{algorithm03}
 is  consistent with that defined by  \eqref{aux0} \eqref{aux1}. We consider    two cases:
1)  $\sigma_{i  , p } = m$. For this case,
 by \eqref{eq111}  \eqref{B20} \eqref{B21} we see $   \tilde{x}_{i,p+1}=  x_{i,p+1} ~ \forall  i:\sigma_{i, p }=m $.
By   $\sigma_{  p+1} = m$ we see  $p+1\in [\tau_m, \tau_{m+1})$,
   and hence  $   \tilde{x}_{i,p+1}=  x_{i,p+1}$ by \eqref{case2}.
So,   the assertion holds  for any $i$ with $  \sigma_{i  , p } = m$.
2)   $\sigma_{i  , p} < m$. From $\sigma_{p+1}=m$ we see  $p+1\in [\tau_m, \tau_{m+1})$,
and hence  by  $\sigma_{i  , p} < m$
from  \eqref{case3} we  see that
 $   \tilde{x} _{i, p+1} $  defined by     \eqref{aux0} \eqref{aux1} equals $x^{*}$.
 By \eqref{A11} \eqref{B21} we derive  $   \tilde{x} _{i, p+1}=x^{*} $,  and hence
    the assertion holds  for any $i$ with $  \sigma_{i  , p } < m$.

   In summary, we have  shown that   $  \tilde{x}_{i, p+1}$ and $\sigma_{p+1}$ generated by  \eqref{algorithm01}-\eqref{algorithm03}
 are consistent with their definitions  \eqref{aux0} \eqref{aux1} \eqref{largest}.  This   completes the proof. \hfill 
 $\blacksquare$

 \section{Proof of Lemma \ref{lemma1}} \label{PL1} 
 For the lemma it suffices to  prove
  \begin{equation}\label{lem1}
  \begin{split}
   & \lim\limits_{T \rightarrow 0} \limsup\limits_{k \rightarrow \infty} \frac{1}{T}
    \parallel \sum_{s=n_k}^{m(n_k,t_k)\wedge ( \tau_{\sigma_{n_k}+1}-1 )} \gamma_s  \tilde{\varepsilon}_{i,s+1}
    I_{[ \parallel \tilde{x}_{i, s}\parallel \leq K ]} \parallel
    \\&~~~~~~~    =0 ~~\forall t_k \in[0,T] \textrm{   for sufficiently large  } K  >0
    \end{split}
  \end{equation}
  along indices $\{ n_k\}$
  whenever  $\{ \tilde{x}_{i, n_k}\}$   converges   for  the  sample path  $\omega$  where   A5  b) holds  for agent $i$.

    We consider the following two cases:

Case 1:     $ \lim\limits_{k \rightarrow \infty} \sigma_{k} =  \sigma < \infty $.
We now show
\begin{equation}\label{infinity}
    \tau_{\sigma+1}=\infty\textrm{  when } \lim\limits_{k \rightarrow \infty} \sigma_{k} =  \sigma .
\end{equation}
Recall that   $\sigma_{k}  $     is defined as
 the largest truncation number among all agents at time $k$,  from $\lim\limits_{k \rightarrow \infty} \sigma_{k} =  \sigma$  we have $
 \sigma_{i, k} \leq \sigma ~~ \forall k \geq 0 ~~ \forall i \in \mathcal{V}.$
 From here  by the definition of $\tau_{i,m}$  it follows that
 $\tau_{i, \sigma+1} = \textrm{inf}  \{k: \sigma_{i,k}= \sigma+1 \} = \infty ~~ \forall i \in \mathcal{V},$
and hence   $\tau_{\sigma+1}=\infty$.  Thus, \eqref{infinity} holds.

 From    \eqref{gap11} we have  $\tilde{\tau}_{i,\sigma}  \leq  \tau_{\sigma}+BD   $,
                       and  hence  by  \eqref{aux1} \eqref{infinity}
 \begin{equation}\label{111}
 \tilde{x}_{i,k} =x_{i,k}, ~  \tilde{\varepsilon}_{i,k+1}=\varepsilon_{i,k+1}  ~~\forall k \geq   \tau_{ \sigma} +BD.
 \end{equation}
 So,  $\parallel \sum_{s= k}^{m( k,t )\wedge ( \tau_{\sigma_{ k}+1}-1 )} \gamma_s  \tilde{\varepsilon}_{i,s+1}
    I_{[ \parallel \tilde{x}_{i, s}\parallel \leq K ]} \parallel \\ =\parallel \sum_{s=k}^{m( k,t ) } \gamma_s   \varepsilon_{i,s+1}    I_{[ \parallel x_{i, s}\parallel \leq K ]} \parallel$   for any $t>0$ and any  sufficiently large $k$.
    Then  we conclude  \eqref{lem1} by A5  b).

  Case 2:   $ \lim\limits_{k \rightarrow \infty} \sigma_{k} =   \infty $.
  We    prove    \eqref{lem1}  separately for  the following three cases:

 i)  $\tilde{\tau}_{i, \sigma_{ n_p } } \leq n_{p } $. For this case,
   $[ n_p,\tau_{\sigma_{n_p }+1}) \subset [ \tilde{\tau}_{i, \sigma_{n_p }},\tau_{\sigma_{n_p }+1}) ,$     and hence from  \eqref{aux1}  we derive
 \begin{equation}\label{C1}
   \tilde{x}_{i,s}= x_{i,s}, ~  \tilde{\varepsilon}_{i,s+1}= \varepsilon_{i,s+1}  ~  \forall s: n_p   \leq s < \tau_{\sigma_{n_p }+1}.
 \end{equation}
Thus, for sufficiently large $K$ and  any $   t_p \in[0,T]$
 \begin{equation}\label{s1}
 \begin{array}{lll}
    &   \parallel \sum_{s=n_p }^{m(n_p ,t_p)\wedge ( \tau_{\sigma_{n_p}+1}-1 )} \gamma_s  \tilde{\varepsilon}_{i, s+1}
    I_{[ \parallel \tilde{x}_{i, s}\parallel \leq K ]} \parallel   \\&  =  \parallel
     \sum_{s=  n_p }^{m(n_p ,t_p)\wedge ( \tau_{\sigma_{n_p }+1}-1)}
     \gamma_s  \varepsilon_{i,s+1}  I_{[ \parallel x_{i,s}\parallel \leq K ]}  .
     \end{array}
 \end{equation}
By \eqref{C1} we conclude  that $ \{ x_{i,n_p } \} $ is a convergent subsequence.
 Noticing $\sum_{s=  n_p }^{m(n_p ,t_p)\wedge ( \tau_{\sigma_{n_p }+1}-1)}  \gamma_s  \leq
 \sum_{s=  n_p }^{m(n_p ,t_p) }  \gamma_s  \leq t_p \leq T ,$
  we then  conclude \eqref{lem1}    by  \eqref{s1} and  A5 b).

ii)     $\tilde{\tau}_{i, \sigma_{ n_p } } >  n_{p } $    and  $\tilde{\tau}_{i, \sigma_{n_p }}=\tau_{\sigma_{n_p }+1}$.
By the  definitions of $\tau_k$ and $\sigma_k$ we derive $\tau_{\sigma_k} \leq k$, and hence
 $\tau_{\sigma_{n_p}} \leq n_p$.
   Then    $[ n_p,\tau_{\sigma_{n_p }+1}) \subset [ \tau_{\sigma_{n_p}} , \tilde{\tau}_{i, \sigma_{n_p }}) $,
and hence  by \eqref{aux0} we have
     \begin{equation}
    \begin{array}{lll}
     & \tilde{x}_{i,s}= x^*, \quad  \tilde{\varepsilon}_{i,s+1}= -f_i(x^{*}) \quad
 \forall s:  n_p  \leq s <\tau_{\sigma_{n_p }+1}.     \nonumber
    \end{array}
    \end{equation}
From  $\tilde{\tau}_{i, \sigma_{n_p }}=\tau_{\sigma_{n_p }+1}$  by   \eqref{gap11} we see
    $ \tau_{\sigma_{n_p }+1} \leq \tau_{\sigma_{n_p }}+BD  \leq n_p  +BD.  $
    Then for sufficiently large $K$ and  any $   t_p \in[0,T]$
     \begin{equation}
 \begin{array}{lll}
    &   \parallel \sum_{s=n_p }^{m(n_p ,t_p)\wedge ( \tau_{\sigma_{n_p}+1}-1 )} \gamma_s  \tilde{\varepsilon}_{i, s+1}
    I_{[ \parallel \tilde{x}_{i, s}\parallel \leq K ]} \parallel   \\&  \leq
     \sum_{s=  n_p }^{ n_p+BD}
     \gamma_s \parallel  f_i(x^*)    \parallel  \xlongrightarrow [p \rightarrow \infty]{}  0 , \nonumber
     \end{array}
 \end{equation}
 and hence \eqref{lem1}  holds.

 iii)
  $\tilde{\tau}_{i, \sigma_{ n_p } } >  n_{p } $ and   $\tilde{\tau}_{i, \sigma_{n_p }} < \tau_{\sigma_{n_p }+1}$.
  By the definition of   $\tilde{\tau}_{i, \sigma_{ n_p } }$  from   $\tilde{\tau}_{i, \sigma_{n_p }} < \tau_{\sigma_{n_p }+1}$
   it follows that $\tilde{\tau}_{i, \sigma_{ n_p } }=\tau_{i, \sigma_{n_p }}   $, and hence  by \eqref{gap11}
   $ \tau_{i, \sigma_{n_p }}   \leq \tau_{\sigma_{n_p }}+BD$.
   By  noticing $\tau_{\sigma_{n_p}} \leq n_p$ we    conclude  that
    \begin{equation}\label{gap0}
 \tau_{\sigma_{n_p}} \leq  n_p  < \tilde{\tau}_{i, \sigma_{ n_p } }= \tau_{i, \sigma_{n_p }}  \leq n_p +BD  .
    \end{equation}
So,     $[ n_p,\tau_{i, \sigma_{n_p }}) \subset [ \tau_{\sigma_{n_p}} , \tilde{\tau}_{i, \sigma_{n_p }}) $.
  Then from here  and $\tilde{\tau}_{i, \sigma_{ n_p } }= \tau_{i, \sigma_{n_p }} $
   by \eqref{aux0} \eqref{aux1} we derive  \begin{equation}\label{s31}
 \begin{array}{lll}
 & \tilde{x}_{i,s}= x^*, \quad \tilde{\varepsilon}_{i,s+1}= -f_i(x^{*}) \quad
 \forall s: n_p  \leq s <\tau_{i, \sigma_{n_p }}, \\
 & \tilde{x}_{i,s}= x_{i,s},  \quad \tilde{\varepsilon}_{i,s+1}= \varepsilon_{i,s+1}  \quad
  \forall s:\tau_{i, \sigma_{n_p }}  \leq s < \tau_{\sigma_{n_p }+1}. \nonumber
  \end{array}
 \end{equation}
Consequently,  for sufficiently large $K$ and  any $   t_p \in[0,T]$
 \begin{equation}\label{s3}
 \begin{array}{lll}
  &     \parallel \sum_{s=n_p }^{m(n_p ,t_p)\wedge  ( \tau_{\sigma_{n_p }+1}-1)}
  \gamma _s  \tilde{\varepsilon}_{i,s+1}
  I_{[ \parallel \tilde{x}_{i,s}\parallel \leq K  ]}  \parallel
   \\& \leq   \parallel \sum_{s=n_p }^{m(n_p ,t_p)\wedge ( \tau_{\sigma_{n_p }+1}-1)} \gamma _s f_i(x^*)
   I_{[ n_p  \leq s <\tau_{i, \sigma_{n_p }}]} \parallel \\
  & +  \parallel \sum_{s=  \tau_{i, \sigma_{n_p }}     }^{m(n_p ,t_p)\wedge ( \tau_{\sigma_{n_p }+1}-1)}
     \gamma_s  \varepsilon_{i,s+1}  I_{[ \parallel x_{i,s}\parallel \leq K ]}  \| .
 \end{array}
 \end{equation}

Note  that the first term at the right hand of  \eqref{s3} is smaller than
$ \sum_{s=n_p }^{ \tau_{i, \sigma_{n_p }} } \gamma _s \|  f_i(x^*)\|$, which tends to zero as $k  \rightarrow \infty$
by the last inequality   in  \eqref{gap0} and $ \lim\limits_{k \rightarrow \infty} \gamma_{k} =  0.$
 The truncation number   for  agent $i$ at time
   $\tau_{i, \sigma_{n_p }}$ is $\sigma_{n_p }$,  while it is smaller than $\sigma_{n_p }$
   at time $\tau_{i, \sigma_{n_p }}-1$  since   $\tau_{i, \sigma_{n_p }}$ is the smallest time
   when the truncation number of $i$ has reached $\sigma_{n_p }$. Consequently, by  Remark \ref{r1} we have
     $x_{i, \tau_{i, \sigma_{n_p }}}=x^*$, and hence    $\{x_{i, \tau_{i, \sigma_{n_p }}} \}_{p \geq 1}$ is a
convergent subsequence.   Noticing $$\sum_{s=  \tau_{i, \sigma_{n_p }}     }^{m(n_p ,t_p)\wedge ( \tau_{\sigma_{n_p }+1}-1)}
     \gamma_s  \leq \sum_{s=  n_p }^{m(n_p ,t_p) }  \gamma_s  \leq t_p  ,$$
  we conclude \eqref{lem1}    by  \eqref{s3} and  A5 b).

Since one of  i),  ii), iii) must take place for the case $ \lim\limits_{k \rightarrow \infty} \sigma_{k} =   \infty $,
we thus have proved    \eqref{lem1}  in Case 2.

 Combining Case 1 and Case 2  we derive  \eqref{lem1}. $\blacksquare$
 
   \section{ Proof  of  Lemma \ref{lemma2}} \label{PL2}
 Let us  consider  a fixed  $\omega $ where A5 b) holds.

Let  $C > \|  \bar{X}\| $.  There exists an integer  $k_C >0 $ such that
   \begin{equation}\label{subsequence}
 \|  \widetilde{X}_{n_k} \| \leq C  \quad \forall k \geq k_C.
   \end{equation}

 By  Lemma \ref{lemma1},   there exist a constant $T_1>0$ and a positive integer  $k_0 \geq k_C$ such that
    \begin{equation}\label{lem23}
    \begin{split}
&   \Big \|  \sum_{s=n_k}^{m(n_k,t_k)\wedge ( \tau_{\sigma_{n_k}+1}-1)} \gamma_s  \tilde{\varepsilon}_{    s+1}
      I_{[\parallel \widetilde{X}_{ s}\parallel \leq  K  ]} \Big \|  \leq  T_0
    \\&  \forall t_k \in[0,T_0] \quad  \forall T_0  \in[0, T_1 ]~~\forall k \geq k_0
    \end{split}
   \end{equation}
for sufficiently large  $K >0$.  Define
     \begin{align}
   & M_0'\triangleq  1+C( c\rho +2),  \label{def21} \\
     &H_1\triangleq  \max_{X} \{ \parallel F(X) \parallel :   \parallel X \parallel \leq M_0'+1
     +C  \},  \label{def22} \\
  &c_1\triangleq H_1+3 +  \frac{c ( \rho+1)  }{1-\rho},~~\mbox{and} \quad c_2\triangleq \frac{H_1+1}{\sqrt{N}},    \label{def23}
\end{align}
  where  $c$ and  $\rho$ are given by \eqref{graphmatrix}. Select  $T >0$  such that
 \begin{equation}\label{def24}
     0< T  \leq T_1\textrm{  and  }c_1 T  <1.
 \end{equation}

For  any $ k \geq k_0$ and any  $ T_k \in [0,T]$,   define
   \begin{equation}\label{lem24}
    \begin{split}
  & s_k\triangleq  \sup \{ s \geq n_k: \parallel \widetilde{X}_j- \widetilde{X}_{n_k} \parallel\\&
   \leq c_1T_k +M_0' \quad \forall j: n_k \leq j \leq s \}.
    \end{split}
    \end{equation}
Then from  \eqref{subsequence} and  \eqref{def24} it follows that
    \begin{equation}\label{lem22}
     \begin{split}
   \parallel \widetilde{X}_s \parallel
    \leq   M_0'+1+C\quad   \forall s: n_k \leq s \leq s_k .
   \end{split}
   \end{equation}

 We intend to prove  $ s_k >  m(n_k, T_k).$
 Assume the converse that for sufficiently large $ k  \geq k_0$ and any $ T_k \in [0,T]$
 \begin{equation}\label{inverse2}
    s_k \leq m(n_k, T_k).
 \end{equation}
  We first  show that  there exists a positive integer  $k_1 > k_0$
   such that for sufficiently large $ k  \geq k_1$
 \begin{equation}\label{notruncation}
      s_k < \tau_{\sigma_{n_k}+1}~~ \forall k \geq  k_1 ~~ \forall T_k \in [0,T].
 \end{equation}
   We prove \eqref{notruncation}  for the  two  alternative   cases:  $ \lim\limits_{k \rightarrow \infty } \sigma_k = \infty$
  and $ \lim\limits_{k \rightarrow \infty } \sigma_k =\sigma <\infty $.

 i)     $ \lim\limits_{k \rightarrow \infty } \sigma_k = \infty$.
 Since $\{M_k\}$ is a sequence of positive numbers increasingly  diverging to infinity,
   there exists a positive  integer  $k_1 > k_0$
   such that    $M_{ \sigma_{n_k} }  >M_0'+1+C$ for all $k \geq k_1$.
 Hence, from \eqref{lem22} we know     $ s_k < \tau_{\sigma_{n_k}+1}  $.

  ii)  $ \lim\limits_{k \rightarrow \infty } \sigma_k =\sigma <\infty $.
  For this  case there exists a positive  integer  $k_1 > k_0$
      such that  $ \sigma_{n_k} =\sigma$ for all    $k \geq k_1 $, and hence
      $\tau_{\sigma_{n_k}+1}=\tau_{\sigma+1}=\infty$ by \eqref{infinity}. Then  $m(n_k,T_k)  < \tau_{\sigma_{n_k}+1}$,
       and hence  by \eqref{inverse2} we derive \eqref{notruncation}.

 By  \eqref{def24},   we see that $ T_k \in [0,T] \subset [0, T_1]$. Then from   \eqref{lem23} it follows that
  for  sufficiently large $ k \geq k_1$ and $K >0$:    \begin{equation}\label{lem208}
  \begin{split}
  &\left \| \sum_{s=n_k}^{m(n_k,t_k )\wedge ( \tau_{\sigma_{n_k}+1}-1)} \gamma_s  \tilde{\varepsilon}_{   s+1}
      I_{[\parallel \tilde{X}_{  s}\parallel \leq  K  ]}  \right\| \leq  T_k
 \\&\qquad~~\qquad~~\qquad~~ \forall t_k  \in[0,T_k]~~\forall  T_k \in [0,T].
   \end{split}
   \end{equation}  Set  $t_k =\sum_{m=n_k}^s \gamma_m $ for some $s \in [ n_k,s_k].$ Then
 $\sum_{m=n_k}^s \gamma_m \leq \sum_{m=n_k}^{s_k} \gamma_m  \leq T_k   $ by \eqref{inverse2}.
Since  $m(n_k,t_k) =s$,  from \eqref{notruncation} we have that $m(n_k,t_k )\wedge ( \tau_{\sigma_{n_k}+1}-1) =s.$
 Then by setting    $ K   \triangleq M_0'+1+C$,
  from   \eqref{lem22}  \eqref{lem208}  it follows that
 for sufficiently large   $   k \geq k_1$ and any $T_k \in [0,T]$:
 \begin{equation}\label{err}
 \begin{split}
 & \left\|  \sum_{m=n_k}^{s} \gamma_m  \tilde{\varepsilon}_{  m+1}  \right \|   \leq   T_k
  \quad \forall s: n_k \leq s \leq s_k.
      \end{split}
    \end{equation}

Let us consider the  following algorithm starting from  $n_k$ without truncation
  \begin{equation}\label{lem26}
Z_{m+1}= (W(m) \otimes \mathbf{I}_l)    Z_m + \gamma_{m } (F(Z_m)+\tilde{\varepsilon}_{m+1}),
~ Z_{n_k}= \tilde{X}_{n_k}.
  \end{equation}
By   \eqref{notruncation}  we know that
  \eqref{compactform}  holds for $m=n_k, \cdots,  s_k-1$ for   $   \forall k \geq k_1 ~ \forall T_k \in [0,T]$.
  Then  we have the following:
  \begin{equation}\label{aux21}
      Z_m =   \widetilde{X}_m  ~~ \forall m: n_k \leq m \leq s_k  .
  \end{equation}
Hence  by \eqref{def22} \eqref{lem22}, we know that for $   \forall k \geq k_1 ~ \forall T_k \in [0,T]$
\begin{equation}\label{upper21}
  \parallel   F(Z_m)   \parallel \leq H_1 ~~ \forall m: n_k \leq m \leq s_k .
\end{equation}
   Set  $z_k = \frac{\mathbf{1}^T \otimes \mathbf{I}_l}{N} Z_k$.
Multiplying both sides of \eqref{lem26}  from left  with $ \frac{1}{N} ( \mathbf{1}^T \otimes \mathbf{I}_l) $,
by  $ \mathbf{1}^TW(m) =\mathbf{1}^T$  and \eqref{kr}, we derive
 $$  z_{s+1}=z_s + \frac{\mathbf{1}^T \otimes \mathbf{I}_l}{N } \gamma_s (F(Z_s)+\tilde{\varepsilon}_{s+1}),$$
 and hence
  \begin{equation}\label{lem27}
  \begin{array}{lll}
 &  \parallel  z_{s+1}  - z_{n_k} \parallel = \parallel \frac{\mathbf{1}^T \otimes \mathbf{I}_l}{N }  \sum_{m=n_k}^s \gamma_{m } (F(Z_m)+\tilde{\varepsilon}_{m+1}) \parallel \\
   &\leq \frac{1}{ \sqrt{N}}  \big(   \sum_{m=n_k}^s \gamma_{m }\parallel   F(Z_m)   \parallel +
    \parallel \sum_{m=n_k}^s   \gamma_{m }  \tilde{\varepsilon}_{m+1}  \parallel  \big). \nonumber
   \end{array}
  \end{equation}
Then  by   \eqref{inverse2}, \eqref{err}, and   \eqref{upper21},
we conclude  that for  sufficiently large $ k \geq k_1$ and any $ T_k \in [0,T]$
   \begin{equation}\label{lem28}
  \begin{split}
  \| z_{s+1}  - z_{n_k} \|  \leq \frac{H_1+ 1}{ \sqrt{N} }  \sum_{m=n_k}^s \gamma_{m }  =c_2  T_k
 ~~ \forall s: n_k \leq s \leq s_k.
   \end{split}
  \end{equation}
Denote by $Z_{\bot, s}\triangleq  D_{\bot} Z_{s}$ the disagreement vector of $Z_{s}$.
By multiplying both sides of  \eqref{lem26}  from left with $ D_{\bot}$,  we derive
  \begin{equation*}\label{dis1}
   \begin{split}
  & Z_{\bot, s+1}=  D_{\bot} (W(s) \otimes I_l)Z_{  s }+ \gamma_{s }  D_{\bot}     (F(Z_s)+\tilde{\varepsilon}_{s+1}),
  \end{split}
   \end{equation*}
 and inductively
   \begin{equation}\label{dis21}
   \begin{array}{lll}
  & Z_{\bot, s+1}= \Psi(s, n_k)  Z_{n_k}+ \\
  &   \sum_{m=n_k} ^s
  \gamma_{m } \Psi(s-1, m) D_{\bot} (F(Z_m)+\tilde{\varepsilon}_{m+1}) \quad \forall s \geq n_k. \nonumber
  \end{array}
   \end{equation} Thus, by   \eqref{power23}   and   \eqref{power24} we  have that
   \begin{small}
 \begin{equation}\label{dis22}
   \begin{array}{lll}
  & Z_{\bot, s+1}   =  [(\Phi(s, n_k) - \frac{1}{N} \mathbf{1} \mathbf{1}^T) \otimes I_l]    Z_{n_k} \\&+ \sum_{m=n_k} ^s
  \gamma_{m }   [(\Phi(s-1, m) - \frac{1}{N} \mathbf{1} \mathbf{1}^T) \otimes I_l]  F(Z_m) \\
  & +\sum_{m=n_k} ^s \gamma_{m }    [( \Phi(s-1, m) - \frac{1}{N} \mathbf{1} \mathbf{1}^T) \otimes I_l]    \tilde{\varepsilon}_{m+1}.
  \end{array}
   \end{equation}
   \end{small}
By \eqref{subsequence} and  \eqref{aux21},  we have that $\| Z_{n_k} \| \leq C$. Thus,
 from  \eqref{graphmatrix} and  \eqref{upper21} it follows that
 \begin{small}
   \begin{equation}\label{dis3}
   \begin{array}{lll}
  & \parallel Z_{\bot, s+1} \parallel \leq C c\rho^{s +1-n_k} + \sum_{m=n_k} ^s
  \gamma_{m }  H_1 c\rho^{s -m} +   \\
  &    \parallel     \sum_{m=n_k} ^s \gamma_{m }    [( \Phi(s-1, m) - \frac{1}{N} \mathbf{1} \mathbf{1}^T) \otimes \mathbf{I}_l]    \tilde{\varepsilon}_{m+1} \parallel .
  \end{array}
   \end{equation}
\end{small}

Define $\Gamma_n \triangleq \sum_{m=1} ^n  \gamma_m    \tilde{\varepsilon}_{m+1}.$
Then by \eqref{err},   we obtain  that
  $\parallel \Gamma_s- \Gamma_{n_k-1} \parallel  \leq   T_k  \quad \forall s: n_k \leq s \leq s_k$.
Note that  \begin{equation} \label{errr3D}
   \begin{array}{lll}
 &   \sum_{m=n_k} ^s
  \gamma_{m }   (\Phi(s-1, m)  \otimes I_l) \tilde{\varepsilon}_{m+1}   \\ &=  \sum_{m=n_k} ^s     (\Phi(s-1, m)  \otimes I_l) (\Gamma_{m}-\Gamma_{m-1})
   \\ &=  \sum_{m=n_k} ^s     (\Phi(s-1, m)  \otimes I_l) (\Gamma_{m}-\Gamma_{n_k-1})
   \\& -\sum_{m=n_k} ^s     (\Phi(s-1, m)  \otimes I_l) (\Gamma_{m-1}-\Gamma_{n_k-1}) . \nonumber
 \end{array}
   \end{equation}
Summing  by parts,  by \eqref{graphmatrix}  we have
 \begin{equation} \label{errr3D}
   \begin{array}{lll}
 &  \left \|  \sum\limits_{m=n_k} ^s
  \gamma_{m }   (\Phi(s-1, m)  \otimes I_l) \tilde{\varepsilon}_{m+1} \right\|\leq \left \| \Gamma_s-  \Gamma_{n_k-1}\right \|
   + \\& \sum\limits_{m=n_k} ^{s-1}     \left \|  \Phi(s-1, m) - \Phi(s-1, m+1)    \right \|
\left \|  \Gamma_{m} - \Gamma_{n_k-1} \right \| \\
  & \leq    T_k +    \sum\limits_{m=n_k} ^{s-1} ( c \rho^{s-m-1}+c \rho^{s-m})  T_k
 \\&  \leq T_k + \frac{c ( \rho+1)  }{1-\rho}  T_k ~~ \forall s: n_k \leq s \leq s_k. \nonumber
 \end{array}
   \end{equation}
This incorporating with \eqref{err}  yields  that for
 sufficiently large $ k \geq k_1$ and any $  T_k \in [0,T]$,
 \begin{small}
\begin{equation}\label{D11f}
\begin{split}
& \| \sum_{m=n_k} ^s \gamma_{m }    [( \Phi(s-1, m) - \frac{1}{N} \mathbf{1} \mathbf{1}^T) \otimes I_l]    \tilde{\varepsilon}_{m+1} \|
\\& \leq (2+\frac{c ( \rho+1)  }{1-\rho})T_k  \quad \forall s: n_k \leq s \leq s_k.
\end{split}
\end{equation}
\end{small}

By noticing $\sum_{m=n_k} ^s  \gamma_{m }   \rho^{s -m}   \leq \frac{ 1 }{1-\rho} \sup_{ m \geq n_k  } \gamma_{m }
 ,$ from   $ \lim\limits_{k \rightarrow \infty} \gamma_{k} =  0 $, \eqref{dis3}, and  \eqref{D11f} it    follows
 that  for sufficiently large $ k \geq k_1$ and any $  T_k \in [0,T]$:
    \begin{equation}\label{est4}
   \begin{array}{lll}
  & \parallel Z_{\bot, s+1} \parallel  \leq  Cc\rho + \frac{ c H_1 }{1-\rho} \sup_{ m \geq n_k  } \gamma_{m }
 + (2 +  \frac{c ( \rho+1)  }{1-\rho} ) T_k
   \\& \leq Cc\rho + 1
 + (2 +  \frac{c ( \rho+1)  }{1-\rho} ) T_k  ~~ \forall s: n_k \leq s \leq s_k.
  \end{array}
   \end{equation}
By noticing
 $Z_s=Z_{\bot,s}+( \mathbf{1} \otimes I_l ) z_s ,$ we derive
  \begin{equation*}\label{D12}
   \begin{array}{lll}
   & \parallel Z_{s+1}  - Z_{n_k} \parallel \\&=
    \parallel ( \mathbf{1} \otimes \mathbf{I}_l)z _{s+1}+ Z_{\bot, s+1}  - Z_{\bot, n_k} -
    ( \mathbf{1} \otimes \mathbf{I}_l ) z_{n_k} \parallel \\
      & \leq  \parallel  Z_{\bot, s+1} \parallel + \parallel  Z_{\bot, n_k} \parallel + \sqrt{N}\parallel
    z_{s+1}- z_{n_k}  \parallel.
  \end{array}
   \end{equation*}
Since  $\parallel  Z_{\bot, n_k} \parallel \leq 2\parallel  Z_{  n_k} \parallel=2C $,
 from \eqref{lem28} and \eqref{est4} it follows   that for sufficiently large $k \geq k_1$  and  any $  T_k \in [0,T]$
 \begin{small}
    \begin{equation}\label{result1}
   \begin{array}{lll}
   & \parallel Z_{s+1}  - Z_{n_k} \parallel  \leq C c\rho  + 1
 + (2+  \frac{c ( \rho+1)  }{1-\rho} ) T_k   +2C + \\&   (H_1+ 1) T_k
    \leq C( c\rho + 2 )+1 + (3+H_1+  \frac{c ( \rho+1)  }{1-\rho} )
    T_k \\& =M_0'+c_1 T_k  \quad \forall s: n_k \leq s \leq s_k.
  \end{array}
   \end{equation}
   \end{small}
Therefore,  from \eqref{def24} and  $\| Z_{n_k} \| \leq C$  we know that  for sufficiently large $k \geq k_1$ and  any $ T_k \in [0,T]$
 \begin{equation}\label{D13}
    \parallel Z_{s_k+1} \parallel \leq \parallel  Z_{n_k} \parallel  +     M_0'+c_1 T_k \leq
M_0'+1+C. \nonumber
 \end{equation}
Rewrite  \eqref{algorithm01} in the  compact form as follows
 $$  \widehat{X}_{s_k+1}= [W(s_k) \otimes \mathbf{I}_l]
   \widetilde{X}_{s_k} + \gamma_{s_k } (F( \widetilde{X}_{s_k})+\tilde{\varepsilon}_{s_k+1}) ,$$
   where $ \widehat{X}_k \deq  col \{ \hat{x}_{1,k}, \cdots, \hat{x}_{N,k} \}.$
   Then    $\widehat{X}_{s_k+1}=Z_{s_k+1}$  by \eqref{lem26}  \eqref{aux21}, and hence from  \eqref{D13}  it follows that
 \begin{equation}\label{bound}
    \parallel  \widehat{X}_{s_k+1} \parallel   \leq  M_0'+1+C  .
 \end{equation}

  We now show
  \begin{equation}\label{appentr}
      \widetilde{X}_{s_k+1}= \widehat{X}_{s_k+1}\textrm{ and }   s_k+1 < \tau_{\sigma_{n_k}+1}
  \end{equation}
 for sufficiently large  $k \geq k_1$ and any $  T_k \in [0,T]$.
   We   consider  the following two cases:  $ \lim\limits_{k \rightarrow \infty } \sigma_k = \infty$
  and $ \lim\limits_{k \rightarrow \infty } \sigma_k =\sigma <\infty $.

i)      $ \lim\limits_{k \rightarrow \infty } \sigma_k = \infty$.
By noting   $M_{ \sigma_{n_k} }  >M_0'+1+C~~ \forall k \geq k_1$,  from \eqref{bound},  \eqref{algorithm02} and
 \eqref{algorithm03} we get  $\widetilde{X}_{s_k+1}=\widehat{X}_{s_k+1}$ and $\sigma_{ s_k+1}=\sigma_{ s_k }.$
Hence  $s_k+1 < \tau_{\sigma_{n_k}+1}$ by \eqref{notruncation}.

ii)      $ \lim\limits_{k \rightarrow \infty } \sigma_k =\sigma< \infty$. Since
$ ~\tau_{\sigma_{n_k}+1}= \infty   ~\forall k \geq k_1$,
 by  \eqref{inverse2}  we see $      s_k+1 < \tau_{\sigma_{n_k}+1} $ . Then by $\sigma_{n_k}=\sigma  ~\forall k \geq k_1$
 we conclude  $\sigma_{s_k+1}=\sigma_{s_k}=\sigma $.
Hence    by \eqref{algorithm02} we derive    $\widetilde{X}_{s_k+1}=    \widehat{X}_{s_k+1}$.
Thus, \eqref{appentr} holds.

From \eqref{appentr}, we know  that  \eqref{compactform}  holds for    $m=s_k$
for sufficiently large  $k \geq k_1$ and any $  T_k \in [0,T]$.
By   $\widehat{X}_{s_k+1}=Z_{s_k+1}$ and \eqref{appentr},  we obtain   $\widetilde{X}_{s_k+1}=Z_{s_k+1}$.
Hence from  \eqref{result1} and  $Z_{n_k}= \widetilde{X}_{n_k}$  it follows that
 for sufficiently large  $k \geq k_1$ and any $  T_k \in [0,T]$
$$\parallel  \widetilde{X}_{s_k+1}
 -  \widetilde{X}_{n_k} \parallel \leq M_0'+c_1 T_k.$$
 This  contradicts  with  the definition of $s_k$ given  by \eqref{lem24}.
Thus, \eqref{inverse2} does not hold.
Consequently, $s_k > m(n_k,T_k)$, and hence  by the definition of   $s_k$ given  in \eqref{lem24}   we  derive \eqref{res21}.

Since $\{ \widetilde{X}_s: n_k \leq s \leq m(n_k,T_k)\}$ are  bounded,
similar to proving   \eqref{notruncation} it can be shown that $m(n_k,T_k) +1   < \tau_{\sigma_{n_k}+1} $.
As a result,     \eqref{compactform}
holds for $m=n_k,\cdots, m(n_k,T_k)$. Similar to    \eqref{lem28},  we have the following
assertion for  sufficiently large $k  $ and any $ T_k \in [0,T]$:
   \begin{equation}
  \begin{split}
 &  \parallel  \bar{x}_{m+1}  - \bar{x}_{n_k} \parallel   \leq  c_2T_k
   \quad \forall m: n_k \leq m \leq m(n_k,T_k)  \nonumber.
   \end{split}
  \end{equation} Hence, \eqref{res22}   holds.
   The proof of  Lemma \ref{lemma2} is completed.
    \hfill $\blacksquare$

    \section{Proof  of Lemma \ref{lemnoise}} \label{PLN}

    Since   $ \lim\limits_{k \rightarrow \infty}\widetilde{X}_{n_k}  = \bar{X}$,  we then have
 $\lim\limits_{k \rightarrow \infty} \bar{x}_{n_k} =  \bar{x}=\frac{\mathbf{1}^T \otimes \mathbf{I}_l}{N} \bar{X}$.
  Lemma \ref{lemma2} ensures that  there exists a $T\in (0,1) $
such that   for sufficiently large  $k$,
$m(n_k,T ) < \tau_{\sigma_{n_k}+1}$ and $\{\widetilde{X}_s: n_k \leq s \leq m(n_k,T ) +1\}$ are bounded.  So,   for  any $ T_k \in[0,T]$ and  any  sufficiently large  $K$
\begin{small}
 \begin{align*}
\left \| \sum_{s=n_k}^{m(n_k,T_k)\wedge ( \tau_{\sigma_{n_k}+1}-1 )} \gamma_s  \tilde{\varepsilon}_{s+1}
    I_{[ \parallel \tilde{X}_ s\parallel \leq K ]} \right \|   = \left \|  \sum_{s=n_k}^{m(n_k,T_k) } \gamma_s  \tilde{\varepsilon}_{ s+1}   \right \| .
\end{align*}
\end{small}
Therefore,  by  Lemma \ref{lemma1} we derive
 \begin{equation}\label{error20}
    \lim_{T \rightarrow 0} \limsup_{k \rightarrow \infty} \frac{1}{T}
    \left \|  \sum_{s=n_k}^{m(n_k,T_k) } \gamma_s  \tilde{\varepsilon}_{ s+1}  \right \|  =0 \quad \forall T_k \in[0,T]. \nonumber
 \end{equation}
So, for   \eqref{err2}  it suffices to show
  \begin{equation}\label{error23}
    \begin{split}
    &   \lim_{T \rightarrow 0} \limsup_{k \rightarrow \infty} \frac{1}{T}
   \left \| \sum_{s=n_k}^{m(n_k,T_k) } \gamma_s   e_{s+1}   \right \|   =0 \quad \forall T_k \in[ 0,T  ].
  \end{split}
    \end{equation}

 Similar to  \eqref{est4},
   we can show  that there exist positive constants  $c_3 ,  c_4   ,c_5  $  such that for sufficiently large  $k$
 \begin{equation}\label{lem31}
  \begin{array}{lll}
      \parallel \tilde{X}_{\bot, s+1} \parallel & \leq  c_3 \rho^{s+1-n_k} +
 c_4  \sup_{ m\geq n_k } \gamma_{m }     \\& + c_5 T \quad \forall s: n_k \leq s \leq  m(n_k,T ).
      \end{array}
   \end{equation}
Since $ 0 < \rho<1$, there exists  a positive integer $m' $ such that $\rho^{m'}<T.$
Then  $\sum_{m=n_k}^{n_k+m'} \gamma_m \xlongrightarrow [k \rightarrow \infty]{} 0 $
 by   $\lim\limits_{ k \rightarrow \infty} \gamma_k =0   $.
Thus,      $ n_k+m'  < m(n_k,T)$  for sufficiently large   $k  $.
Then by  \eqref{lem31},   we know  that   for sufficiently large   $k  $
 \begin{equation}\label{bd1}
 \parallel \tilde{X}_{\bot, s+1} \parallel \leq   o(1)   + (c_3+c_5) T
  \quad \forall s:  n_k+m' \leq s \leq  m(n_k,T ),
   \end{equation}
where    $o(1) \rightarrow 0$ as $k \rightarrow \infty$.

 Since $ \bar{x}_{n_k} \xlongrightarrow [k \rightarrow \infty]{}  \bar{x}$,
  from \eqref{res22} and \eqref{bd1} it follows that  for sufficiently large   $k  $ and
  $\forall s:  n_k+m' \leq s \leq  m(n_k,T )$
   \begin{equation}\label{bound12}
    \begin{split}
      & \parallel   \bar{x}_s  -\bar{x} \parallel     \leq     \parallel   \bar{x}_s - \bar{x}_{n_k} \parallel+ \parallel \bar{x}_{n_k} -\bar{x} \parallel=o(1) + \delta(T), \\
    & \parallel \tilde{x}_{i,s} -\bar{x} \parallel  \leq    \parallel \tilde{x}_{i,s} - \bar{x}_s \parallel+
    \parallel   \bar{x}_s - \bar{x}  \parallel =o(1) +\delta(T)  , \nonumber
  \end{split}
    \end{equation}
   where  $ \delta(T) \rightarrow 0$ as $T \rightarrow 0$. By continuity  of  $f_i(\cdot)$ we derive
  $ \parallel f_i(\tilde{x}_{i,s}) -f_i(\bar{x}_s) \parallel   \leq \parallel f_i(\tilde{x}_{i,s}) -f_i(\bar{x} ) \parallel  +\parallel f_i(\bar{x}_s)-f_i(\bar{x})   \parallel  =o(1) +  \delta(T).$
Consequently,
$$\parallel e_{i,s+1} \parallel =   \parallel f_i(\tilde{x}_{i,s}) -f_i(\bar{x}_s) \parallel /N =o(1) +  \delta(T).$$
Then for sufficiently large   $k  $
     \begin{equation}\label{D00}
         \parallel  e_{s+1}  \parallel =o(1) +  \delta(T)
         \quad \forall s:  n_k+m' \leq s \leq  m(n_k,T ).
     \end{equation}

By the boundedness of   $\left \{\widetilde{X}_s: n_k \leq s \leq m(n_k,T ) \right \}$  and the continuity
    of   $f_i(\cdot)$,  we know that there exists a constant $c_e >0$  such that
 $ \left \| e_{s +1} \right \|\leq c_e $. Then from  \eqref{D00}  we derive
 \begin{small}
    \begin{equation}
    \begin{split}
    &    \left \| \sum_{s=n_k}^{m(n_k,T_k) } \gamma_s   e_{s+1}   \right \|    \leq    \sum_{s=n_k}^{ n_k+m' }\gamma_s c_e  +
  ( o(1) +  \delta(T) )  \sum_{s=  n_k+m'+1}^{m(n_k,T_k) } \gamma_s
   \\& \leq         c_e m' \sup_{s \geq n_k} \gamma_s +    T(  \delta(T) +o(1)  ~~\forall T_k \in [0,T]  \nonumber
  \end{split}
    \end{equation}
    \end{small}
  for sufficiently large   $k  $. Then by   $\lim\limits_{ k \rightarrow \infty} \gamma_k =0   $, we derive
      \begin{equation}\label{error22}
    \begin{split}
    &    \limsup_{k \rightarrow \infty } \frac{1}{T} \parallel \sum_{s=n_k}^{m(n_k,T_k) } \gamma_s   e_{s+1}  \parallel
  =  \delta(T)  \quad \forall T_k \in[ 0,T  ]. \nonumber
  \end{split}
    \end{equation}
 Letting $T\rightarrow 0$, we  derive  \eqref{error23}, and hence \eqref{err2} holds.

   The proof of  Lemma \ref{lemnoise} is completed.
    \hfill $\blacksquare$

\section{  Proof of Lemma \ref{lemma4}}\label{PL4}

  i)   Assume $
    \lim\limits_{k \rightarrow \infty } \sigma_k =  \infty.$ Then  there exists a positive integer $ n_k $
  such that  for   any $  k\geq 1$, $ \sigma_{n_k}=k $ and  $\sigma_{n_k-1}=k-1$.
Hence $n_k-1 \in [\tau_{k-1},\tau_k)$ by  $ \sigma_{n_k}=k $. Thus, from \eqref{sub2}
  it follows that  $\tilde{x}_{i,n_k} =x^*~\forall i \in \mathcal{V}$.
Consequently,    $ \widetilde{X}_{n_k}= (\mathbf{1} \otimes I_l) x^{*}$  and hence
$\{ \widetilde{X}_{n_k}\}$  is a convergent subsequence  with     $\bar{x}_{n_k}= x^{*}$.

  Since    $\{M_k \}$  is a    sequence of positive numbers increasingly diverging to infinity,
   there exists  a positive integer  $k_0   $  such that
   \begin{equation}\label{no2}
    M_{ k } \geq  2\sqrt{N}c_0+ 2+M_1' \quad \forall k \geq k_0,
   \end{equation}
where $c_0$ is given in  A2  c) and
 \begin{equation}\label{F2}
    M_1'=2+  ( 2\sqrt{N}c_0 +2)( c \rho +2).
 \end{equation}

In what follows, we show that under the converse assumption    $\{\bar{x}_{n_k}\}$   starting from  $ x^{*}$
 crosses the sphere with  $ \| x \|=c_0$  infinitely many times.
   Define
     \begin{align}
    &     m_k \triangleq \inf\{ s >  n_k: \parallel \widetilde{X}_s \parallel  \geq 2\sqrt{N}c_0+2+M_1' \},  \label{notion11} \\
     &   l_k \triangleq \sup\{ s < m_k: \parallel \widetilde{X}_s \parallel  \leq 2\sqrt{N}c_0 +2 \}.  \label{notion12}
\end{align}
Noticing  $ \|  \widetilde{X}_{n_k} \| =\sqrt{N} \|  x^{*}\|$
    and $\|  x^{*}\|<  c_0$, we derive $ \| \widetilde{X}_{n_k} \| < \sqrt{N} c_0.$
Hence from    \eqref{notion11}   \eqref{notion12}  it is seen  that  $n_k <  l_k  < m_k .$
By the definition of  $l_k$, we know that  $\{ \widetilde{X}_{l_k} \}$ is bounded. Then there
exists a convergent subsequence,  denoted still by    $\{ \widetilde{X}_{l_k}\}.$
 By   $\bar{X} $ denoting  the limiting point of  $\widetilde{X}_{l_k} $, from
  $\| \widetilde{X}_{l_k} \|  \leq 2\sqrt{N}c_0 +2 $
  it follows that  $ \| \bar{X} \| \leq  2\sqrt{N}c_0+2$.

 By  Lemma   \ref{lemma2}, there exist constants     $M_0'>0$ defined by \eqref{def21} with  $C= 2\sqrt{N}c_0+2 $,
  $c_1>0  $ and $c_2>0$ defined by \eqref{def23},    $ 0< T<1$ with    $c_1 T \leq 1$ such that
  \begin{equation}
       \parallel  \widetilde{X}_{m+1} -\widetilde{X}_{l_k}  \parallel \leq c_1 T +M_0'
   \quad \forall m: l_k \leq m \leq m(l_k,T)  \nonumber
  \end{equation}
for sufficiently large   $k \geq k_0$. Then   for sufficiently large   $k \geq k_0$
  \begin{equation}\label{res71}
  \begin{split}
    & \parallel  \widetilde{X}_{ m+1}   \parallel     \leq  \parallel  \widetilde{X}_{l_k}  \parallel + c_1T +M_0' \\
    &  \leq 2\sqrt{N}c_0+ 2+1+1+(  2\sqrt{N}c_0+2) (c\rho+2)
  \\& =2\sqrt{N}c_0+2+M_1'  \quad \forall m: l_k \leq m \leq m(l_k,T).
  \end{split}
  \end{equation}
 Then  $m(l_k,T) \leq n_{k+1}$ for sufficiently large $k \geq k_0 $ by \eqref{no2}.

From \eqref{res71}  by the definition of  $m_k$ defined in \eqref{notion11}, we conclude  $ m(l_k,T)+1 \leq m_k$
    for sufficiently large   $k \geq k_0.$ Then by
     \eqref{notion11} \eqref{notion12}, we know that for sufficiently large   $k \geq k_0$
 \begin{equation}\label{llm3}
 \begin{split}
    &  2\sqrt{N}c_0 +2 < \parallel \widetilde{X}_{m+1} \parallel \leq 2\sqrt{N}c_0+2+M_1'
\\&~~~~~~~~~~~~~~~~~~~~~~\quad \forall m: l_k \leq m \leq m(l_k,T) .
\end{split}
 \end{equation}
     Since $0< \rho<1$, there exists a positive integer  $m_0$  such that  $4c\rho^{m_0} <1$.
   Then  $\sum_{m=l_k}^{l_k+m_0} \gamma_m \to 0 $ by $ \gamma_k\to 0$, and hence    $ l_k+m_0  < m(l_k,T) < n_{k+1} $   for sufficiently large   $k \geq k_0$.
So, from  \eqref{llm3} it is seen that  for sufficiently large   $k \geq k_0$
      \begin{equation}\label{F4}
        \parallel  \widetilde{X}_{ l_k+m_0} \parallel >   2\sqrt{N}c_0+2.
      \end{equation}

Noticing that  $ \{\widetilde{X}_{m+1} : l_k \leq m \leq m(l_k,T)\}$  are bounded,
similarly  to  \eqref{est4} we know  that  for sufficiently large   $k \geq k_0$
    \begin{equation}  \label{llem2}
   \begin{split}
  & \parallel  \widetilde{X}_{\bot, m+1} \parallel \leq ( 2\sqrt{N}c_0+2)c\rho^{m+1-l_k} +
  \frac{ c H_1 }{1-\rho} \sup_{ m \geq l_k  } \gamma_{m }    \\
    & + 2  T +  \frac{c ( \rho+1)  }{1-\rho} T
    \quad \forall m: l_k \leq m \leq m(l_k,T). \nonumber
  \end{split}
   \end{equation}
From here, by $c_1T<1 $ and $\gamma_k \xlongrightarrow [k \rightarrow \infty]{} 0$ it follows that
     \begin{equation}    \label{llem3}
   \begin{split}
  & \parallel  \widetilde{X}_{\bot, l_k+m_0} \parallel \leq  2 c\rho^{m_0}(  \sqrt{N}c_0 +1) + \frac{1}{2}+ c_1
   T \\ & \leq \frac{1}{2}(  \sqrt{N}c_0 +1) + \frac{1}{2}+1 =\frac{\sqrt{N}c_0}{2}+2
  \end{split}
   \end{equation}
    for sufficiently large   $k \geq k_0$.  By noticing
    $  (\mathbf{1} \otimes \mathbf{I}_l )\bar{x}_{l_k+m_0}= \widetilde{X}_{ l_k+m_0}- \widetilde{X}_{\bot,l_k+m_0} ,$
from   \eqref{F4} and   \eqref{llem3}  we conclude that
     \begin{equation}
   \begin{split}
&   \sqrt{N} \parallel \bar{x}_{l_k+m_0} \parallel=  \parallel \widetilde{X}_{ l_k+m_0}- \widetilde{X}_{\bot, l_k+m_0} \parallel \\&
\geq  \parallel \widetilde{X}_{ l_k+m_0} \parallel -  \parallel\widetilde{X}_{\bot, l_k+m_0} \parallel > \frac{3}{2}\sqrt{N}c_0.
  \end{split}
   \end{equation}
Therefore,   $\parallel \bar{x}_{l_k+m_0} \parallel >c_0.$

   Thus, we have shown that for sufficiently large   $k \geq k_0$, starting from
    $ x^{*},$   $\{\bar{x}_{n_k}\}$  crosses  the sphere with  $ \| x \| =c_0$ before the time $n_{k+1}$.
So,   $\{\bar{x}_{n_k}\}$   starting from $ x^{*}$  crosses   the sphere with   $\| x \| =c_0$  infinitely many times.

 ii) Since   $v(J)$  is nowhere dense, there exists a nonempty interval  $[\delta_1, \delta_2 ]
    \in (v(x^*),\textrm{inf }_{\parallel x \parallel =c_0} v(x))$  with  $d([\delta_1, \delta_2 ],  v(J))>0 $.
   Assume the converse $\lim\limits_{k \rightarrow \infty } \sigma_k =  \infty.$ By i)  $\{v(\bar{x}_k)\}$ crosses the interval  $[\delta_1, \delta_2 ]$
      infinitely many times while $ \{\widetilde{X}_{n_k}\}$ converges. This contradicts Lemma \ref{lemma3}.
   Therefore, the inverse assumption  $    \lim\limits_{k \rightarrow \infty } \sigma_k =  \infty$   is impossible,
    and  hence \eqref{lem4} holds. The proof is completed.
  \hfill $\blacksquare$

\end{document}